\documentclass[review]{siamltex}
\usepackage[english]{babel}
\usepackage{amsmath,amssymb,amsfonts}
\usepackage{color}
\usepackage{graphicx}
\usepackage{lipsum}

\newcommand{\qed}{\hfill$\Box$}

\newcommand{\R}{\mathbb R}
\newcommand{\N}{\mathbb N}

\newcommand{\oma}{\Omega}
\newcommand{\xinto}{\int_\oma}

\newcommand{\texinto}{\int_0^T\!\!\int_\oma}
\newcommand{\cs}{\mathcal{S}}
\newcommand{\cl}{\mathcal{L}}
\newcommand{\pt}{\partial_t}
\newcommand{\vp}{\varphi}
\newcommand{\bs}{\overline{s}}
\newcommand{\by}{\overline{y}}

\renewcommand{\le}{\leqslant}

\renewcommand{\ge}{\geqslant}

\begin{document}

\title{A new type of identification problems: \\optimizing the 
fractional order \\ in a nonlocal 
evolution equation}

\author{J\"urgen Sprekels\thanks{Department of Mathematics, Humboldt-Universi\-t\"at zu Berlin,
Unter den Linden 6, 10099 Berlin, Germany, and Weierstrass Institute
for Applied Analysis and
Stochastics, Mohrenstrasse 39, 10117 Berlin, Germany.}
\and
Enrico Valdinoci\thanks{Weierstrass Institute
for Applied Analysis and
Stochastics, Mohrenstrasse 39, 10117 Berlin, Germany,
and Dipartimento di Matematica Federigo Enriques,
Universit\`a degli Studi di Milano,
Via Cesare Saldini 50, 20133 Milano, Italy.}
}
\maketitle

\begin{abstract}
In this paper, we consider a rather general linear evolution equation of fractional 
type, namely a diffusion type problem in which the diffusion operator is the $s$th power of a positive definite operator having a discrete spectrum in $\R^+$. We prove existence, uniqueness and differentiability properties with respect
to the fractional parameter $s$. These results are then employed to derive 
existence as well as
first-order necessary and second-order sufficient optimality conditions
for a minimization problem, which is inspired by considerations in mathematical biology. 

In this problem,
the fractional parameter $s$ serves as the ``control parameter'' that  
needs to be chosen in such a way as to minimize a given cost functional. This problem constitutes a new class
of identification problems: while usually in identification problems the type of the
differential operator is prescribed and one or several of its coefficient functions
need to be identified, in the present case one has  to determine the type 
of the differential operator itself. 

This problem exhibits the inherent analytical 
difficulty that
with changing fractional parameter $s$ also the domain of definition, and thus 
the underlying function space,  of the
fractional operator changes.    
\end{abstract}



\section{Introduction}

\noindent
Let $\Omega\subset \R^n$ be a given open domain and, with a given
$T>0$, $Q:=\oma\times (0,T)$. We consider in $\oma$ the evolution of 
a fractional diffusion process governed by the 
$s-$power of a positive definite operator ${\mathcal{L}}$. In this paper,
we study, for a given $L\in(0,+\infty)\cup \{+\infty\}$,  the following identification problem for fractional
evolutionary systems:

\vspace{3mm}\noindent
{{\bf (IP)} \,\,Minimize the cost function} 
\begin{align}
\label{cost}
J(y,s):=
\frac{1}{2}
\texinto \big|y(x,t)-y_Q(x,t)\big|^2\,dx\,dt
\,+\,\vp(s)\end{align}
with~$s$ in the interval $(0,L)$, subject to the fractional evolution problem
\begin{align}
\label{ss1}
& {\pt y\,+\,\cl^s y\,=\,f \quad\mbox{in }\,Q,}\\[1mm]
\label{ss2}
& {\hspace*{5mm} y(\cdot,0)\,=\,y_0 \quad\mbox{in }\,\oma.}
\end{align}

\vspace{2mm}
\noindent
In this connection, $y_Q\in L^2(Q)$ is a given target function, and 
$\vp\in C^2(0,L)$ is a nonnegative penalty function satisfying   
\begin{align}
\label{phi}
\lim_{s\searrow 0}\vp(s)=+\infty=\lim_{s\nearrow L}\vp(s).
\end{align}
Examples of penalty functions which fulfill \eqref{phi} are
\begin{equation}\label{EXA:PHI}
\begin{split}
&\vp(s)=\frac{1}{s\,(L-s)}\qquad{\mbox{ for $s\in(0,L)$, if $L\ne+\infty$}}\\
{\mbox{and }}
&\vp(s)=\frac{e^s}{s}\qquad{\mbox{ for $s\in(0,L)$,
if $L=+\infty$}}.
\end{split}
\end{equation}
The properties of the right-hand side $f$ and of the initial datum $y_0$
will be specified later.

{Problem {\bf (IP)} defines a class of identification problems which,
to the authors' best knowledge, has never been studied before. Indeed, while there
exists a vast literature on the identification of coefficient functions or of
right-hand sides in parabolic and hyperbolic evolution equations (which cannot be cited 
here), there are only but a few contributions to the control theory of 
fractional operators of diffusion type. In this connection, we refer the reader to the
recent papers \cite{AO1}, \cite{AO2}, \cite{AO3} and \cite{Bors}. However, in these
works the fractional operator was fixed and given a priori.
% Also in the recent identification paper \cite{CL}, a coefficient function of a fractional operator
%of prescribed type (i.\,e., $s$ is given) was sought. 
In contrast to these
papers, in our case the type of the fractional order operator itself, which is defined
by the parameter $s$, is to be determined.} 

The fact that the fractional order parameter $s$ is the ``control variable''
in our problem entails a mathematical difficulty, namely, that with changing $s$ 
also the domain of $\cl^s$ changes. As a consequence, in the functional analytic
framework also the underlying solution space changes with $s$. From this,
mathematical difficulties have to be expected. For instance, simple compactness
arguments are likely not to work if existence is to be proved. 
In order to overcome this difficulty, we present in Section 4
(see the compactness result of Lemma 6) an argument which is based on Tikhonov's
compactness theorem.

Another feature of the problem {\bf (IP)} is the following:
if we want to establish necessary and sufficient optimality conditions, then
we have to derive differentiability properties of the control-to-state
($s\mapsto y$) mapping. A major part of this work is devoted to
this analysis.   

In this paper, the fractional power of the diffusive operator is seen as an ``optimization
parameter''. This type of problems has natural applications. For instance,
a biological motivation is the following: in the study of the diffusion of biological species
(see, e.g., \cite{RO, F12, P13, MV} and the references therein)
 there is experimental evidence
(see~\cite{W96, H10}) that many predatory species
follow ``fractional'' diffusion patterns instead of classical
ones: roughly speaking, for instance, suitably long excursions
may lead to a more successful  hunting strategy.
In this framework, optimizing over the fractional parameter~$s$
reflects into optimizing over the ``average excursion''
in the hunting procedure, which plays a crucial role
for the survival and the evolution of a biological population
(and, indeed, different species in nature adopt different
fractional diffusive behaviors).

In this connection\footnote{As a technical remark, we point out that,
strictly speaking, in view of their probabilistic
and statistical interpretations, many of the experiments 
available in the literature are often more
closely related to fractional operators of integrodifferential
type rather than to fractional operators of spectral
type, and these two notions are, in general, not
the same (see e.g.~\cite{TWO}), although they coincide,
for instance, on the torus, and are
under reasonable assumptions asymptotic to each other
in large domains
(see e.g. Theorem~1 in~\cite{NAZA3} for precise estimates).
Of course, the problem considered in this paper 
does not aim to be
exhaustive, and other types of operators and
cost functions may be studied as well,
and, in fact, in concrete situations
different ``case by case''
analytic and phenomenological considerations
may be needed to produce detailed
models which are as accurate as possible
for ``real life'' applications.}, the solution $y$ to the state system
\eqref{ss1}, \eqref{ss2} can be thought of as the spatial
density of the predators (where the birth and death rates of the 
population are not taken into account here, but rather
its capability of adapting to the environmental situation).
In this sense, the minimization of $\,J\,$ is related
to finding the ``optimal''
distribution for the population (for instance, in terms
of the availability of resources, possibility of using favorable environments,
distributions of possible preys, favorable conditions for reproduction, etc.).
Differently from the existing literature, this optimization
is obtained here by changing the nonlocal diffusion parameter~$s$, where,
roughly speaking, a small~$s$ corresponds to a not very dynamic
population and a large~$s$ to a rather mobile one.

The growth condition \eqref{phi} has to be understood against this
biological background: in nature, neither a complete immobility of
the individuals (i.\,e., the choice $s=0$) nor an extremely fast
diffusion (observe that even the extreme case $s=L=+\infty$  
is allowed in our setting) are 
likely to guarantee the survival of the species. In this connection, we
may interpret the target function $y_Q$ as, e.\,g., the spatial distribution 
of the prey. To adapt their strategy, the predators must know these seasonal
distributions a priori; however, this is often the case from long standing
experience. We also remark that in nature the prey species in turn adapt their
behavior 
to the strategy of the predators; it would thus be more realistic
to consider a predator-prey system with two (possibly different) values
of $s$. Such an analysis, however, goes beyond the scope of this work
in which we confine ourselves to the simplest possible situation.

The remainder of the paper is organized as follows: in the following section,
we formulate the functional analytic framework of our problem and prove  the basic well-posedness
results for the state system \eqref{ss1}, \eqref{ss2}, as well as 
its differentiability properties with respect to the
parameter $s$. Afterwards, in Section 3, we study the problem {\bf (IP)} and establish the 
first-order necessary and the second-order sufficient 
conditions of optimality. Some elementary explicit
examples are also provided, in order to show the influence
of the boundary data and of the target distribution on the
optimal exponent.

The final section then brings an existence result whose proof employs  
a compactness result (established in Lemma 6), which is
based on Tikhonov's compactness theorem.

\section{Functional analytic setting and results for the solution operator}

\noindent
The mathematical setting in which we work is the following: we consider an
open and bounded  domain~$\Omega\subset\R^n$
and a differential operator~${\mathcal{L}}$
acting on functions mapping~$\Omega$ into~$\R$, together
with appropriate boundary conditions.
We generally assume that there exists a complete
orthonormal system (i.\,e., an orthonormal basis) $\{e_j\}_{j\in\N}$
of~$L^2(\Omega)$ having the property that each~$e_j$ lies
in a suitable subspace~${\mathcal{D}}$ of~$L^2(\Omega)$, 
and such that~$e_j$ is an eigenfunction of~${\mathcal{L}}$
with corresponding 
eigenvalue~$\lambda_j\in\R$, for any~${j\in\N}$
(notice that in this way
the boundary conditions of the differential
operator~${\mathcal{L}}$ can be encoded in
the functional space~${\mathcal{D}}$).
In this setting, we may write, {for any $j\in\N$,} 
$$ \mathcal{L} e_j =\lambda_j e_j {\mbox{ in }}\Omega,\quad
e_j\in{\mathcal{D}}.$$
We also generally assume that \begin{equation*}{\mbox{$\lambda_j\ge0$
for any~$j\in\N$.}}\end{equation*}
The prototype of operator~${\mathcal{L}}$ that we have
in mind is,
of course, (minus) the Laplacian in a bounded
and smooth domain~$\Omega$ (possibly
in the distributional sense), together with either Dirichlet or Neumann
homogeneous boundary conditions (in these cases, for smooth domains,
one can take, respectively,
either~${\mathcal{D}}:=H^2(\Omega)\cap
H^1_0(\Omega)$ or~${\mathcal{D}}:=
H^2(\Omega)$).

For any~$v$, $w\in L^2(\Omega)$, 
we consider the scalar product
$$ \langle v, w\rangle := \int_\Omega v(x)\,w(x)\,dx.$$
In this way,
we can write any function~$v\in L^2(\Omega)$ {in the form}
$$ v =\sum_{j\in\N} \langle v, e_j\rangle\, e_j,$$
where the equality is indented in the~$L^2(\Omega)$-sense,
and, if
$$ v\in {\mathcal{H}}^1 := \Big\{
v\in L^2(\Omega):\,\,
\{\lambda_j \,\langle v, e_j\rangle\}_{j\in\N}\in\ell^2
\Big\}$$
then
$$ {\mathcal{L}} v = \sum_{j\in\N} \lambda_j \,
\langle v, e_j\rangle\,e_j.$$
For any~$s>0$, we define the $s$-power of the operator~${\mathcal{L}}$
in the following way. First, we consider the space
\begin{equation}\label{HD}
{\mathcal{H}}^s := \left\{v\in L^2(\Omega):\,\,{\|v\|_{\mathcal{H}^s}<+\infty}
\right\},\end{equation}
where we use the notation
\begin{equation}\label{Su-s-2s}
\| v \|_{{\mathcal{H}}^s}:=
\Bigl(\sum_{j\in\N} \lambda_j^{2s} \,\big|\langle v, e_j\rangle\big|^2\Bigr)^{1/2}.
\end{equation}
Notice that the notation of the space ${\mathcal{H}}^s$
has been chosen in such a way
that ${{\mathcal{H}}^s}$, for $s=1$, reduces to the space
${{\mathcal{H}}^1}$ that was introduced above.
This notation is reminiscent of, but different from, the notation for fractional
Sobolev spaces
(roughly speaking, $s=1$ in our notation forces
the Fourier coefficients to be in $\ell^2$
weighted by one power of the eigenvalues; in the case
of second order operators this would correspond to Sobolev spaces
of order two, rather than one, and this difference in the notation
is the main reason for which we chose
to use calligraphic fonts for our functional spaces).

{We then set, for any~$v\in{\mathcal{H}}^s$,}
\begin{equation}
\label{Ls}
{\mathcal{L}}^s v:= \sum_{j\in\N} \lambda_j^s \,
\langle v, e_j\rangle\,e_j.
\end{equation}

We are ready now to define our notion of a solution to the state system: given~$y_0\in L^2(\Omega)$
and~$f:\Omega\times[0,T]\to\R$ such that~$f(\cdot,t)\in L^2(\Omega)$
for every~$t\in[0,T]$, 
we say that~$\,y:\Omega\times[0,T]\to\R\,$
is a solution to the state system \eqref{ss1}, \eqref{ss2}, if and only if the following
conditions are satisfied:
\begin{eqnarray}
\label{SL1}&&\mbox{$y(\cdot, t)\in{\mathcal{H}}^s$ \,for any \,$t\in(0,T]$,}\\
\label{OV}&& \lim_{t\searrow0} \,\langle y(\cdot,t),\,e_j\rangle=
\langle y_0,\,e_j\rangle \quad\mbox{for all }\,j\in\N,\\
\label{SL2}&&\mbox{for every $j\in\N$, the mapping $(0,T)\ni t\mapsto
\langle y(\cdot,t),\,e_j\rangle$ is } \\
&&\nonumber\mbox{absolutely continuous,}\\
\label{SL3}&&{\mbox{and it holds\,\,
$\partial_t \langle y(\cdot,t),\,e_j\rangle
+\lambda_j^s \langle y(\cdot,t),\,e_j\rangle
=\langle f(\cdot,t),\,e_j\rangle,
$}} \\
&&\nonumber\mbox{for every $j\in\N$ and almost every $t\in (0,T)$.}
\end{eqnarray}
We remark that conditions~\eqref{SL1}, \eqref{OV},
\eqref{SL2} and~\eqref{SL3} are precisely the
functional analytic translations of
the functional identity in~\eqref{ss1}, \eqref{ss2}.

\vspace{5mm}
{We begin our analysis with a result that establishes 
existence, uniqueness and regularity  of the solution to
the state system \eqref{ss1}, \eqref{ss2}. }

\begin{theorem}\label{TH}
\,\,\,Suppose that $f:\Omega\times[0,T]\to\R$ satisfies~$f(\cdot,t)\in L^2(\Omega)$, 
for every~$t\in[0,T]$, as well as
\begin{equation}\label{L2-l2}
\sum_{j\in\N} f_j^2\,<\,+\infty, \quad\mbox{where }
\,f_j:=\sup_{\theta\in(0,T)} \big| \langle f(\cdot,\theta),e_j\rangle\big|\,.
\end{equation}
Then the following holds true:

\vspace{2mm}\noindent
{\bf (i)} \,\,\,\, If $\,y_0\in L^2(\Omega),$ then there exists for every $s>0$
 a unique solution \,$y(s):=y\,$ 
to the state system {\rm \eqref{ss1}, \eqref{ss2}} that 
fulfills the conditions {\rm \eqref{SL1}--\eqref{SL3} {\em and belongs to} $L^2(Q)$}.
Moreover, with the control-to-state operator $\cs:s\mapsto y(s)$, we have the 
explicit representation
\begin{equation}\label{sol1}
\cs (s)(x,t)=y(s)(x,t) = \sum_{j\in \N} y_j(t,s)\,e_j(x) \quad\mbox{a.\,e. in }\,Q,
\end{equation}
where, for $j\in\N$ and $t\in [0,T]$, we have set
\begin{equation}\label{sol2}
y_j(t,s):= \langle y_0,\,e_j\rangle\,e^{-\lambda_j^s t}
+\int_0^t 
\langle f(\cdot,\tau),\,e_j\rangle \,e^{\lambda_j^s (\tau-t)}\,d\tau.\end{equation}

\vspace{2mm}\noindent
{\bf (ii)} \,\,\,If $\,y_0\in \mathcal{H}^{s/2}$, then 
\begin{equation}
\label{regu}
y(s)\in H^1(0,T;L^2(\Omega))\cap L^\infty(0,T;{\mathcal{H}}^{s/2})\cap L^2(0,T;{\mathcal{H}}^s)\,,
\end{equation}
and
\begin{equation}
\label{dty}
\partial_t y(s)\,=\,\sum_{j\in\N} \partial_t  y_j(\cdot,s)\,e_j\,.
\end{equation}
Moreover,  we have the estimate
\begin{align}\label{la:sec2}
&\| \partial_t y(s)\|_{L^2(Q)}^2\,+\,\|y(s)\|^2_{L^\infty(0,T;{\mathcal{H}^{s/2}})}\,
+\,\| y(s)\|_{L^2(0,T;{\mathcal{H}}^s)}^2\,\\[1mm]
&\nonumber \le\,T\sum_{j\in\N}\sup_{\theta\in (0,T)}|\langle f(\cdot,\theta), e_j\rangle|^2\,+\,
\|y_0\|_{{\mathcal{H}}^{s/2}}^2.
\end{align}
\end{theorem}

\vspace{5mm}\noindent
{\bf Remark:} \,\,\,We point out that
formula~\eqref{sol2} is of classical flavor and related
to Duhamel's Superposition Principle. In our setting, this kind
of explicit representation is 
an auxiliary tool used to prove the regularity
estimates with respect to the fractional parameter $s$ that will be needed
later in this paper.\medskip

\vspace{5mm} \noindent
{\sc Proof of Theorem~\ref{TH}:} \,\,\,{\bf (i)}: 
\,\,\,We first prove that the series defined in \eqref{sol1} represents a function in $L^2(Q)$.
To this end, we show that \,$\{\sum_{j=1}^n y_j(\cdot,s)\,e_j\}_{n\in\N}$\, forms a Cauchy 
sequence in $L^2(Q)$. Indeed, we have, for every $n,p\in\N$, the identity
\begin{align}
\label{Cauchy1}
&\left\|\mbox{$\sum_{j=1}^{n+p}y_j(\cdot,s)\,e_j\,-\,\sum_{j=1}^n y_j(\cdot,s)$}\,e_j\right\|^2_{L^2(Q)}
\\
\nonumber &=\,\int_0^T\left\|\mbox{$\sum_{j=n+1}^{n+p} y_j(t,s)\,e_j$}\right\|^2_{L^2(\Omega)} dt
\,=\,\int_0^T\sum_{j=n+1}^{n+p}|y_j(t,s)|^2 dt\,.
\end{align} 
Now, for any $\tau\in(0,t)$, we have that
$ e^{\lambda_j^s (\tau-t)} \le 1$,
since $\lambda_j\ge0$. Accordingly,
\begin{eqnarray*}
&& \left|
\int_0^t
\langle f(\cdot,\tau),\,e_j\rangle \,e^{\lambda_j^s (\tau-t)}\,d\tau\right|
\le \int_0^t \big|
\langle f(\cdot,\tau),\,e_j\rangle\big| \,e^{\lambda_j^s (\tau-t)}\,d\tau
\\ &&\qquad \le \int_0^t \big|
\langle f(\cdot,\tau),\,e_j\rangle\big| \,d\tau
\le T\sup_{\theta\in (0,T)}
|\langle f(\cdot,\theta),e_j\rangle|.\end{eqnarray*}
Thus, it follows from \eqref{sol2} that for every $j\in\N$ and $t\in [0,T]$ it holds
\begin{align*}
|y_j(t,s)|\,\le\,|\langle y_0,e_j\rangle| \,+\,T\sup_{\theta\in (0,T)}
|\langle f(\cdot,\theta),e_j\rangle|\,.
\end{align*}
Since $y_0\in L^2(\Omega)$, we have $\,\sum_{j\in\N}|\langle y_0,e_j\rangle|^2=\|y_0\|_{L^2(\Omega)}^2$, 
and it readily follows from \eqref{L2-l2} that the sequence 
$\{ \sum_{j=1}^n \int_0^T |y_j(t,s)|^2 dt\}_{n\in\N}$\, is a Cauchy sequence in $\R$, which
proves the claim.

\vspace{3mm}
Next,
we observe that
\begin{eqnarray}
\label{estif}
\sup_{\theta\in(0,T)} \| f(\cdot,\theta)\|_{L^2(\Omega)}^2 &=&
\sup_{\theta\in(0,T)} \sum_{j\in\N} \big| \langle f(\cdot,\theta),e_j\rangle\big|^2
\\
\nonumber &\le&\sum_{j\in\N} \sup_{\theta\in(0,T)}\big| \langle f(\cdot,\theta),e_j\rangle\big|^2,
\end{eqnarray}
which is finite, thanks to~\eqref{L2-l2}.
Consequently,
\begin{equation}\label{L2-finite}
\begin{split}
& \int_0^T\| f(\cdot,t)\|^2_{L^2(\Omega)}\,dt < +\infty\\
{\mbox{and }}\quad&
\int_0^T\| f(\cdot,t)\|_{L^2(\Omega)}\,dt < +\infty.
\end{split}\end{equation}

Now, we prove the asserted existence result
by showing that the function \,$y(s)$, which is explicitly defined 
by \eqref{sol1}, \eqref{sol2} in the statement of the theorem, fulfills for every $s>0$ all of the
conditions \eqref{SL1}--\eqref{SL3}.
To this end, let $s>0$ be fixed. We
set, for $j\in\N$ and $t\in [0,T]$,
\begin{align}
\label{V_W_DE}
 v_j(t,s):=\langle y_0,\,e_j\rangle\,e^{-\lambda_j^s t},\quad
 w_j(t,s):=
\displaystyle\int_0^t \langle f(\cdot,\tau),\,e_j\rangle 
\,e^{\lambda_j^s (\tau-t)}\,d\tau.
\end{align}

\noindent
Since $y(s)\in L^2(Q)$, we conclude from \eqref{sol1} and \eqref{sol2} that for every 
$j\in\N$ and $t\in [0,T]$ it holds that
\begin{align}\label{U2}
\langle y(s)(\cdot,t),e_j\rangle\,&
=\,\lim_{n\to\infty} \sum_{k=1}^n\langle y_k(t,s)\,e_k,e_j\rangle\\[1mm]
\nonumber &=\,y_j(t,s)\,=\,v_j(t,s)+w_j(t,s).
\end{align}
Moreover, for any~$t\in(0,T]$, we set
$$ \kappa(t) := \sup_{r\ge 0}\left(r e^{-rt}\right).$$
Notice that~$\kappa(t)<+\infty$ for any~$t\in(0,T]$, and 
\begin{eqnarray*}
\lambda_j^s \,| v_j(t,s)|&\le&\lambda_j^s 
\big|\langle y_0,\,e_j\rangle\big|\,e^{-\lambda_j^s t}\,\le\, \kappa(t)\, 
\big|\langle y_0,\,e_j\rangle\big|\,.
\end{eqnarray*}
Since~$y_0\in L^2(\Omega)$, we therefore have
\begin{equation}\label{U3}
\mbox{$\{\lambda_j^s\, v_j(t,s)\}_{j\in\N}\in\ell^2$, for any~$t\in(0,T]$.}
\end{equation}
In addition, it holds that
\begin{eqnarray*}
\lambda_j^s \,| w_j(t,s)|&\le&
\int_0^t \big|\langle f(\cdot,\tau),\,e_j\rangle\big|\,\lambda_j^s 
\,e^{\lambda_j^s (\tau-t)}\,d\tau\\
&\le& \sup_{\theta\in (0,T)}
\big|\langle f(\cdot,\theta),\,e_j\rangle\big|
\, \int_0^t 
\lambda_j^s
\,e^{\lambda_j^s (\tau-t)}\,d\tau \\
&=& \sup_{\theta\in (0,T)}
\big|\langle f(\cdot,\theta),\,e_j\rangle\big|\,(1-e^{-\lambda_j^s t})
\\ &\le& \sup_{\theta\in (0,T)}
\big|\langle f(\cdot,\theta),\,e_j\rangle\big|,\end{eqnarray*}
and we infer from~\eqref{L2-l2} that also
$\,\{\lambda_j^s\, w_j(t,s)\}_{j\in\N}\in\ell^2\,$, for any~$t\in(0,T]$.
Combining this with \eqref{U2} and~\eqref{U3}, we see that also the sequence \,$
\{\lambda_j^s \,
\langle y(s)(\cdot,t),e_j\rangle\}_{j\in\N}$ belongs to $\ell^2$,
for any~$t\in(0,T]$. 
Thus, by~\eqref{HD} and \eqref{Su-s-2s}, we conclude that $\,y(s)(\cdot,t)\in
{\mathcal{H}}^s\,$
for any~$t\in(0,T]$, and this proves~\eqref{SL1}.

\vspace{2mm}
Next, we point out that~\eqref{OV}
follows directly from~\eqref{sol2},
and thus we focus on the proof of~\eqref{SL2}
and~\eqref{SL3}. To this end, fix~$t\in(0,T)$. If~$|h|>0$ is so small that $t+h\in (0,T)$,
then we observe that
\begin{align}\label{KH:8a8a-1}
&w_j(t+h,s)-w_j(t,s)\\
&=\, e^{-\lambda_j^s(t+h) }
\int_t^{t+h} \langle f(\cdot,\tau),\,e_j\rangle 
\,e^{\lambda_j^s \tau}\,d\tau\nonumber\\
&\quad
+\big( e^{-\lambda_j^s h} -1\big) 
\int_0^{t} \langle f(\cdot,\tau),\,e_j\rangle \,
e^{\lambda_j^s (\tau-t)}
\,d\tau .\nonumber
\end{align}
On the other hand, if we set
$$ g_j(t,s):=
\langle f(\cdot,t),\,e_j\rangle
\,e^{\lambda_j^s t},$$
then we have that
\begin{eqnarray*}
\| g_j(\cdot,s)\|_{L^1(0,T)}
&\le& e^{\lambda_j^s T}\,\int_0^T
\big|\langle f(\cdot,t),\,e_j\rangle\big|\,dt
\\ &\le& 
e^{\lambda_j^s T}\,\int_0^T
\| f(\cdot,t)\|_{L^2(\Omega)}\,dt,
\end{eqnarray*}
which is finite, thanks to~\eqref{L2-finite}.
Hence,  
\begin{equation*}
g_j(\cdot,s)\in L^1(0,T),
\end{equation*} 
and so~$w_j(\cdot,s)$ is absolutely
continuous, and,
by the Lebesgue Differentiation Theorem
(see e.g. \cite{KO} and the references therein),
\begin{eqnarray*}
&& \lim_{h\to0}
\frac1{h}\int_t^{t+h} \langle f(\cdot,\tau),\,e_j\rangle
\,e^{\lambda_j^s \tau}\,d\tau
= \lim_{h\to0} \frac1{h}\int_t^{t+h}g_j(\tau,s)\,d\tau
\\&&\qquad =g_j(t,s)=\langle f(\cdot,t),\,e_j\rangle
\,e^{\lambda_j^s t},\end{eqnarray*}
for almost every $t\in(0,T)$. From this and~\eqref{KH:8a8a-1}, we infer that
$$ \lim_{h\to0}  \frac{
w_j(t+h,s)-w_j(t,s)}{h}
= \langle f(\cdot,t),\,e_j\rangle
-\lambda_j^s
\int_0^{t} \langle f(\cdot,\tau),\,e_j\rangle \,
e^{\lambda_j^s (\tau-t)}
\,d\tau ,$$
for almost every $t\in(0,T)$.
Since also~$v_j(\cdot,s)$ is obviously absolutely continuous, we thus obtain
that $y_j(\cdot,s)$ is absolutely continuous  
and thus differentiable almost everywhere in~$(0,T)$, and we have the identity 
\begin{eqnarray*}
&&\partial_t\langle y(s) (\cdot,t),e_j\rangle\,=\,\partial_t y_j(t,s)\\
&&=\,
-\lambda_j^s\,\langle y_0,\,e_j\rangle\,e^{-\lambda_j^s t}\,+\,
\langle f(\cdot,t),\,e_j\rangle
\,-\,\lambda_j^s
\int_0^{t} \langle f(\cdot,\tau),\,e_j\rangle \,
e^{\lambda_j^s (\tau-t)}
\,d\tau \\[1mm]
&&=\, -\lambda_j^s \,y_j(t,s)+ \langle f(\cdot,t),\,e_j\rangle\\[2mm]
&&=\,-\lambda_j^s \,\langle y(s)(\cdot,t),e_j\rangle\,+\,
\langle f(\cdot,t),\,e_j\rangle, \quad\mbox{for almost every }\,t\in (0,T)\,.
\end{eqnarray*}
This proves~\eqref{SL2} and~\eqref{SL3}.

\vspace{5mm}
As for the uniqueness result, we again fix $s>0$ and assume that there
are two solutions~$y(s), \tilde y(s)\in L^2(Q)$. We put $y^*(s):=
y(s)-\tilde y(s)$, and, adapting the notation
of~\eqref{sol2}, $y^*_j(t,s):=\langle y^*(s)(\cdot,t),e_j\rangle$, for $j\in\N$.
Then, using~\eqref{OV},
\eqref{SL2}, and~\eqref{SL3},
we infer that for every $j\in\N$ the mapping $\,\,t\mapsto y^*_j(t,s)\,\,$
is absolutely continuous in $(0,T)$,  and it satisfies
\begin{equation}\label{SP-01}
\partial_t y^*_j(t,s)
+\lambda_j^s \,y^*_j(t,s)
=0 \quad\mbox{for almost every $\,t\in (0,T)$,}
\end{equation}
as well as
\begin{equation}\label{SP-02}
\lim_{t\searrow0} y^*_j(t,s)=0 .
\end{equation}
Owing to the absolute continuity of $y^*_j(\cdot,s)$, we obtain
(see, e.g., Remark~8 on page~206 of~\cite{brezis}) that~$y^*_j(\cdot,s)\in
W^{1,1}(0,T)$, so that we can use the chain rule 
(see, e.g., Corollary~8.11 in~\cite{brezis}). Thus,
if we define~$\,\zeta_j:=\ln\big(1+(y^*_j(\cdot,s))^2\big)$
and make use of~\eqref{SP-01}, we have that
$$ \partial_t \zeta_j = \frac{2 y^*_j(\cdot,s)\, \partial_t y^*_j(\cdot,s)}
{1+(y^*_j(\cdot,s))^2}
= \frac{-2\lambda_j^s\,  (y^*_j(\cdot,s))^2 }{1+(y^*_j(\cdot,s))^2}\le0
\quad\mbox{a.\,e. in}\,(0,T).$$
%%%% and in the sense  of distributions.
Integrating this relation (see, e.g., Lemma~8.2 in~\cite{brezis}),
we find that, for any~$t_1<t_2\in(0,T)$,
$$ \zeta_j(t_2)\le \zeta_j(t_1).$$
Thus, from~\eqref{SP-02},
$$ \zeta_j(t_2) \le \lim_{t_1\searrow0} 
\zeta_j(t_1) = \lim_{t_1\searrow0}
\ln\big(1+(y^*_j(t_1,s))^2\big) =\ln (1) =0,$$
for any~$t_2\in(0,T)$.
Since also~$\zeta_j\ge0$, we infer that~$\zeta_j$
vanishes identically, and thus also ~$y^*_j(\cdot,s)$.
This proves the uniqueness
claim.

\vspace{5mm}  It remains to show the validity of the claim {\bf (ii)}.
To this end, let again $s>0$ be fixed and assume that $y_0\in {\mathcal{H}}^{s/2}$, which means that
$y_0\in L^2(\Omega)$ 
and 
$$\sum_{j\in\N} \lambda_j^s\,|\langle y_0,e_j\rangle|^2 < +\infty.$$ Now recall
that  \,$\partial_t y_j(t,s)+\lambda_j^s\,y_j(t,s)=\langle f(\cdot,t),e_j\rangle$, for every $j\in\N$
and almost every $t\in (0,T)$. Squaring this equality, we find that
\begin{equation}
|\partial_t y_j(t,s)|^2\,+\,\lambda_j^s\,\frac d{dt} \,|y_j(t,s)|^2 \,+\,\lambda_j^{2s}\,|y_j(t,s)|^2\,=\,
|\langle f(\cdot,t), e_j\rangle|^2\,,
\end{equation} 
and integration over $[0,\tau]$, where $\tau\in [0,T]$, yields that for every $j\in\N$ we have the identity
\begin{align}
&\int_0^\tau |\partial_t y_j(t,s)|^2 dt\,+\,\lambda_j^s\,|y_j(\tau,s)|^2\,+\,\int_0^\tau \lambda_j^{2s}\,|y_j(t,s)|^2 dt
\\
&\nonumber =\,\lambda_j^s\,|\langle y_0,e_j\rangle |^2\,+\,\int_0^\tau| \langle f(\cdot,t), e_j\rangle|^2 dt\,,
\end{align}
whence, for every $n\in \N\cup \{0\}$, $p\in\N$, and $\tau\in [0,T]$,
\begin{align}
\label{Cauchy2}
&\int_0^\tau \sum_{j=n+1}^{n+p}\!\!|\partial_t y_j(t,s)|^2 dt\,+\sum_{j=n+1}^{n+p}\!\!\lambda_j^s\,|y_j(\tau,s)|^2
\,+\int_0^\tau\sum_{j=n+1}^{n+p}\!\! \lambda_j^{2s}\,|y_j(t,s)|^2 dt\\
&\le\,\sum_{j=n+1}^{n+p}\lambda_j^s\,|\langle y_0,e_j\rangle|^2\,+\,T\sum_{j=n+1}^{n+p}
\sup_{\theta\in (0,T)}\,|\langle f(\cdot,\theta),e_j\rangle|^2\,.\nonumber
\end{align}
We remark that we exchanged here summations and integrals: since, up
to now, we are only dealing with a finite summation, this
exchange is valid due to the finite
additivity of the integrals (in particular, we do not need here
any fine result of measure theory).
 
Now, using the same Cauchy criterion argument as in the beginning of the proof of {\bf (i)}, we can therefore infer
that the series 
\begin{align*}
\sum_{j\in\N} \partial_t y_j(\cdot,s)\,e_j, \quad \sum_{j\in\N} \lambda_j^{s/2}
\,y_j(\cdot,s)\,e_j,
\quad\mbox{and }\,\sum_{j\in\N} \lambda_j^s \,y_j(\cdot,s)\,e_j,
\end{align*}
are strongly convergent in the spaces $L^2(Q)$, $L^\infty(0,T;L^2(\Omega))$, and $L^2(Q)$, in this order.
Consequently, we have $\,y(s)\in L^\infty(0,T;{\mathcal{H}}^{s/2})\cap L^2(0,T;\mathcal{H}^s)$,
and also $\partial_t y(s)\in L^2(0,T;L^2(\Omega))$.

We now show that \eqref{dty} holds true, where we denote the limit of series on the right-hand side
by $z$. From the above considerations, we know that, as $n\to\infty$,
$$
\sum_{j=1}^n y_j(\cdot,s)\,e_j\,\to\,y(s), \quad \sum_{j=1}^n \partial_t y_j(\cdot,s)\,e_j\,\to\,z,
\quad\mbox{strongly in }\,L^2(Q)\,.
$$
Hence, there is a subsequence $\,\{n_k\}_{k\in\N}\subset\N\,$ such that, for every test function $\,\phi\in
C^\infty_0(Q)$,
$$
\phi\,\sum_{j=1}^{n_k} y_j(\cdot,s)\,e_j\,\to\,\phi\,y(s),
\quad \phi\,\sum_{j=1}^{n_k} \partial_t y_j(\cdot,s)\,e_j\,\to\,\phi\,z,\quad\mbox{as }\,k\to\infty,
$$
pointwise almost everywhere in $Q$. Using Lebesgue's Dominated Convergence Theorem and Fubini's Theorem
twice, we therefore have the chain of equalities
\begin{align*}
&\int_0^T\!\!\int_\Omega \phi(x,t)\,z(x,t)\,dx\,dt\,=\,\lim_{k\to\infty}\int_\Omega\sum_{j=1}^{n_k} e_j(x)
\int_0^T\!\!\phi(x,t)\,\partial_t y_j(t,s)\,dt\,dx\\
&=\,-\,\lim_{k\to\infty}\int_\Omega\sum_{j=1}^{n_k} e_j(x)\int_0^T\!\!\partial_t\phi(x,t)\,y_j(t,s)\,dt\,dx\\
&=\,-\int_0^T\!\!\int_\Omega\!\partial_t\phi(x,t)\,y(s)(x,t)\,dx\,dt\,,
\end{align*} 
for every $\,\phi\in C^\infty_0(Q)$, that is, we have $\,z=\partial_t y(s)\,$ in the sense of distributions.
Since $z\in L^2(Q)$, it therefore holds $\,y(s)\in H^1(0,T;L^2(\Omega))\,$ with $\,\partial_t y(s)=z$, as claimed. 

Finally,  we obtain the estimate \eqref{la:sec2} from choosing $n=0$ and 
letting $p\to\infty$ in
\eqref{Cauchy2}, which concludes the proof of the assertion. \qed

\vspace{5mm}
Next, we prove an auxiliary result on 
the derivatives of a function of exponential type that
will play an important role in the subsequent analysis. To this end, we define, for fixed 
$\lambda>0$ and~$t>0$,  the real-valued function
\begin{equation}\label{DEF:E}
E_{\lambda,t} (s):= e^{-\lambda^s t} \quad\mbox{for \,$s>0$},
\end{equation}
and denote its first, second, and third derivatives with respect to $s$ by $E_{\lambda,t}' (s)$,
$E_{\lambda,t}'' (s)$, and $E_{\lambda,t}''' (s)$, respectively.
We have the following result.

\vspace*{3mm}
\begin{lemma}\label{DER:1}
There exist constants $\,\widehat{C}_{i}>0$, $0\le i\le 3$, such that, for all
$\lambda>0$, $t\in (0,T]$, and $s>0$, 
\begin{align*}
&\big| E_{\lambda,t} (s)\big|\,\le\,\widehat C_0, \quad 
\left|E_{\lambda,t}'(s)\right|\,\le\,s^{-1}\,\widehat C_1 \,\big(1+|\ln(t)|\big),\\
&\left|E_{\lambda,t}''(s)\right|\,\le\,s^{-2}\,\widehat C_2\,\big( 1+|\ln(t)|^2\big), \quad 
\left|E_{\lambda,t}'''(s)\right|\le\,s^{-3}\,\widehat C_3\,\big(1+|\ln(t)|^3\big)\,.
\end{align*}
\end{lemma}

\vspace{2mm}\noindent
{\sc Proof:} \,\,\,Obviously, we may choose $\widehat C_0=1$, and a simple 
differentiation exercise shows that the 
first three derivatives of $E_{\lambda,t}$ are given by
\begin{align}
&\nonumber E'_{\lambda,t}(s)=-\lambda^s\, t \,e^{-\lambda^s t}\,\ln (\lambda), \quad
E''_{\lambda,t}(s)=\lambda^s\,t\,e^{-\lambda^s t}\left(\lambda^s t-1\right)\,(\ln (\lambda))^2,\\
&\nonumber E'''_{\lambda,t}(s)=\lambda^s\,t\,e^{-\lambda^s t}\left(3\lambda^s t-1-(\lambda^s t)^2\right)
(\ln (\lambda))^3.  
\end{align}

\vspace{1mm}\noindent
Now, observe that
$$ \frac{\ln(\lambda^s t)-\ln (t)}{s} =
\frac{\ln(\lambda^s)+\ln( t)-\ln(t)}{s}=
\ln(\lambda) .$$
Accordingly, we may substitute for $\ln(\lambda)$
in the above identities to obtain that
\begin{align}\label{tbwb-2}
E_{\lambda,t}' (s) \,&= - s^{-1}\lambda^s t\,e^{-\lambda^s t}\,
\big( \ln(\lambda^s t)-\ln(t)\big),
\\
E_{\lambda,t}'' (s) \,&=s^{-2}
\lambda^s t\,e^{-\lambda^s t}\,(\lambda^s t-1)\,
\big( \ln(\lambda^s t)-\ln(t)\big)^2,
\nonumber\\
E_{\lambda,t}'''(s)\,&=s^{-3}\lambda^s t\,e^{-\lambda^s t}
\left(3 \lambda^s t-1-(\lambda^s t)^2\right)
\big(\ln(\lambda^s t)-\ln(t)\big)^3\,.\nonumber 
\end{align}
Thus, we may consider $\,r:= \lambda^s t\,$ as a ``free variable''
in~\eqref{tbwb-2}. Using the fact that 
$$|\ln(r)-\ln(t)|^k\,\le\,2^k\left(|\ln(r)|^k+|\ln(t)|^k\right) \quad\mbox{for }
\,1\le k\le 3,$$
and introducing the finite quantities
\begin{align*}
M_1 &:= \sup_{r>0} \left(r\,e^{-r}\,|\ln (r)|\right),\\
M_2 &:= \sup_{r>0} \left(r\,e^{-r}\right),\\
M_3 &:= \sup_{r>0} \left(r\,e^{-r}\, |r-1|\,4\,|\ln(r)|^2\right),
\end{align*}
\begin{align*}%\\
M_4&:= \sup_{r>0} \left(r\,e^{-r}\, 4\,|r-1|\right),\\
M_5&:= \sup_{r>0} \left(r\,e^{-r}\left|3r-1-r^2\right|\,8\,|\ln(r)|^3\right),\\
M_6&:= \sup_{r>0} \left(r\,e^{-r}\,8\,\left|3r-1-r^2\right|\right),
\end{align*}
we deduce from~\eqref{tbwb-2} the estimates
\begin{align*}
 \big| E_{\lambda,t}' (s)\big|\,&\le\, s^{-1} \big( M_1+M_2\,|\ln (t)|\big),
\\
\big| E_{\lambda,t}'' (s)\big|\,&\le\, s^{-2} \big( M_3+M_4\,|\ln (t)|^2\big),\\
\big|E_{\lambda,t}'''(s)\big|\,&\le\,s^{-3}\big(M_5+M_6\,|\ln(t)|^3\big),
\end{align*}
whence the assertion follows.
\qed

\vspace{5mm}
We are now in the position to derive differentiability properties for the 
control-to-state mapping $\mathcal{S}$. 
As a matter of fact, we will focus on the first and second derivatives,
but derivatives of higher order may be taken into account
with similar methods.
In detail, we have the following result:
\vspace{5mm}

\begin{theorem}\label{II O}
Suppose that
that $f:\Omega\times[0,T]\to\R$ satisfies~$f(\cdot,t)\in L^2(\Omega)$, for every~$t\in[0,T]$, as well
as the condition \eqref{L2-l2}. Moreover, let $y_0\in L^2(\Omega)$. Then the control-to-state
mapping $\mathcal{S}$ is twice Fr\'echet differentiable on $(0,+\infty)$ when viewed as a mapping
from $(0,+\infty)$ into $L^2(Q)$, and for every $\bs\in(0,+\infty)$ the  first and second Fr\'echet
derivatives $D_s\mathcal{S}(\bs)\in\mathcal{L}(\R,L^2(Q))$ and 
$D^2_{ss}\mathcal{S}(\bs)\in\mathcal{L}(\R,\mathcal{L}(\R,L^2(Q)))$ can be identified with the
$L^2(Q)$--functions
\begin{align}
\label{deriS1}
\partial_s y(\bs)&
\,:=\,\sum_{j\in\N}\partial_s y_j(\cdot,\bs)\,e_j,\quad
\partial^2_{ss} y(\bs)\,:=\,\sum_{j\in\N}\partial^2_{ss} y_j(\cdot,\bs)\,e_j\,,
\end{align}
respectively. More precisely, we have, for all $h,k\in\R$,
\begin{equation}
\label{deriS2}
D_s\mathcal{S}(\bs)(h)\,=\,h\,\partial_s y(\bs)\quad\mbox{and }\,
D^2_{ss}\mathcal{S}(\bs)(h)(k)\,=\,h\,k\,\partial^2 _{ss}y(\bs)\,.
\end{equation}
Moreover, there is a constant $\widehat C_4>0$ such that for
all $\bs\in(0,+\infty)$
it holds that
\begin{align}
\label{estideri1}
\left\| D_s\mathcal{S}(\bs) \right\|_{\mathcal{L}(\R,L^2(Q))}&\,=\,\|\partial_s
y(\bs)\|_{L^2(Q)}\,\,\le\,\, \frac{\widehat C_4}{\bs}\,,\\
\label{estideri2}
\left \|D^2_{ss}\mathcal{S}(\bs)\right\|_{\mathcal{L}(\R,\mathcal{L}(\R,L^2(Q)))}&\,=\,
\|\partial^2_{ss} y(\bs)\|_{L^2(Q)}\,\le\,\frac{\widehat C_4}{\bs^2}\,.
\end{align} 
\end{theorem}

\vspace{3mm}\noindent 
{\sc Proof:} \,\,\,Let $\bs\in(0,+\infty)$ be fixed. We first show that 
the functions defined in
\eqref{deriS1} do in fact belong to $L^2(Q)$. To this end, we first note that 
$$
e^{\lambda_j^s (\tau-t)}\,\le\,e^{\lambda_j^s(t-\tau)}\quad \mbox{for }\,0 \le \tau < t, 
$$
and that for $1\le k\le 3$ the functions
\begin{align*}
\phi_k(t):= 1+\,|\ln(t)|^k,\quad \psi_k(t):=\int_0^t\!\!\left(1+|\ln(t-\tau)|^k\right)d\tau,
\quad t\in (0,T],
\end{align*} 
belong to $L^2(0,T)$.

To check this fact, we use the substitution $\theta=-\ln(t)$ and we observe that
\begin{eqnarray*}
&& \int_0^T |\ln(t)|^{k}\,dt\le
\int_0^{T+1} |\ln(t)|^{k}\,dt \\&& \qquad \le \int_0^{1} |\ln(t)|^{k}\,dt +
\int_1^{T+1} |\ln(T+1)|^{k}\,dt \\&& \qquad=
\int_1^{+\infty} \theta^{k}\,e^{-\theta}\,d\theta +
T\,|\ln(T+1)|^{k}\le C(k,T),\end{eqnarray*}
for some $C(k,T)\in(0,+\infty)$. Accordingly,
$$ \int_0^T |\phi_k(t)|^2\,dt
\le 4 \int_0^T \big(1+|\ln(t)|^{2k}\big)\,dt\le 4T +4C(2k,T),$$
hence $\phi_k\in L^2(0,T)$.

Similarly, for any $t\in(0,T]$,
$$ \int_0^t |\ln(t-\tau)|^{k}\,d\tau
= \int_0^t |\ln(\vartheta)|^{k}\,d\vartheta\le
\int_0^T |\ln(\vartheta)|^{k}\,d\vartheta\le C(k,T),$$
and therefore 
$$ |\psi_k(t)|\le T+C(k,T),$$
which gives that $\psi_k\in L^\infty(0,T)\subset L^2(0,T)$, as desired.

Next, we infer from~\eqref{V_W_DE}
and Lemma~\ref{DER:1} that, for every $t\in (0,T]$,
$j\in\N$, and $1\le k\le 3$,
the estimates
\begin{align*}
&\left|\frac{\partial^k}{\partial s^k}\,v_j(t,\bs)\right|\,\le\,
|\langle y_0,e_j\rangle|\left| \frac{d^k}{d s^k}\,E_{\lambda_j,t}(\bs) \right|
\,\le\,\frac{\widehat C_k}{\bs^{\,k}}\,\phi_k(t)\,|\langle y_0,e_j\rangle|\,,\\[1mm]
&\left|\frac{\partial^k}{\partial s^k}\,w_j(t,\bs)\right|\,\le\,\int_0^t
|\langle f(\cdot,\tau), e_j\rangle|\,\left| \frac{d^k}{ds^k}\,E_{\lambda_j,\tau-t}(\bs)\right| 
d\tau \\[2mm]
&\hspace*{25mm}\le\,\widehat C_k\,\bs^{\,-k}\,\psi_k(t)\,\sup_{\theta\in (0,T)}\,|\langle f(\cdot,\theta),
e_j\rangle|\,.\nonumber
\end{align*}

\noindent
Therefore, recalling~\eqref{U2},
we find that, for every $p\in\N$, $n\in\N\cup\{0\}$, and $1\le k\le 2$,  
\begin{align*}
&\left\|\sum_{j=n+1}^{n+p}\frac{\partial^k}{\partial s^k} \,y_j(t,\bs)\,e_j\right\|^2
_{L^2(Q)}
\,\le\,\sum_{j=n+1}^{n+p}\int_0^T\left|\frac{\partial ^k}{\partial s^k}\,y_j(t,\bs)\right|^2 dt\\
&\,\le\,2\sum_{j=n+1}^{n+p}\int_0^T\left|\frac{\partial ^k}{\partial s^k}\,v_j(t,\bs)\right|^2 dt
\,+\,2\sum_{j=n+1}^{n+p}\int_0^T\left|\frac{\partial ^k}{\partial s^k}\,w_j(t,\bs)\right|^2 dt
\nonumber\\
&\nonumber \le\,2\,\widehat C_k^2\,\bs^{\,-2k}\left( \int_0^T\phi_k^2(t)\,dt
\sum_{j=n+1}^{n+p}|\langle y_0,e_j\rangle|^2\right.\\
&\nonumber\hspace*{10mm}\left.+ \int_0^T \psi_k^2(t)\,dt\,\sum_{j=n+1}^{n+p}\,\sup_{\theta\in (0,T)}
\,|\langle f(\cdot,\theta),e_j\rangle|^2\right)\,\,\longrightarrow 0\,, 
\end{align*}
as $n\to\infty$. The Cauchy criterion for series then shows the validity of our claim.
Moreover, taking $n=0$ and letting $p\to\infty$ in the above estimate, we find that  
\eqref{estideri1}
and \eqref{estideri2} are valid provided that \eqref{deriS2} holds true.

\vspace{3mm}
It remains to show the differentiability results. To this end, let $0<|h|<\bs/2$.
Then $\frac 1{\bs-|h|}<\frac 2 {\bs}$, and, invoking Lemma 2
and Taylor's Theorem, we obtain for all $j\in\N$ and $t\in (0,T]$ the estimates
\begin{align*}
&\left|E_{\lambda_j,t}(\bs+h)-E_{\lambda_j,t}(\bs)-h\,E_{\lambda_j,t}'(\bs)\right|
\,=\,\frac 1 2\,h^2\,\left|E_{\lambda_j,t}''(\xi_h)\right|\\[1mm]
&\nonumber\le\,\frac 1 2 \,\widehat C_2\,\xi_h^{-2}\,\phi_2(t)\,h^2
\,\le\, 2 \,\widehat C_2\,\bs^{\,-2}\,\phi_2(t)\,h^2,\\[2mm]
&\left|E_{\lambda_j,t}'(\bs+h)-E_{\lambda_j,t}'(\bs)-h\,E_{\lambda_j,t}''(\bs)\right|
\,=\,\frac 1 2\,h^2\,\left|E_{\lambda_j,t}'''(\eta_h)\right|\\
&\le\,4\,\widehat C_3\,\bs^{\,-3}\,\phi_3(t)\,h^2,\nonumber
\end{align*}
with suitable points $\,\xi_h,\eta_h\in (\bs-|h|,\bs+|h|)$. By the same token,
\begin{align*}
&\int_0^t \left|E_{\lambda_j,\tau-t}(\bs+h)-E_{\lambda_j,\tau-t}(\bs)-
h\,E_{\lambda_j,\tau-t}'(\bs)\right| d\tau\\
&\nonumber\le\,2\,\widehat C_2\,\bs^{\,-2}\int_0^t\phi_2(t-\tau)\,d\tau\,h^2,\\[2mm]
&\int_0^t \left|E_{\lambda_j,\tau-t}'(\bs+h)-E_{\lambda_j,\tau-t}'(\bs)-
h\,E_{\lambda_j,\tau-t}''(\bs)\right| d\tau\\
&\nonumber\le\,4\,\widehat C_3\,\bs^{\,-3}\int_0^t\phi_3(t-\tau)\,d\tau\,h^2\,.
\end{align*}

\noindent
From this, we conclude that with suitable constants $K_i>0$, $1\le i\le 4$, which
depend on $\bs$ but not on $0<|h|<\bs/2$, $j\in\N$, and $t\in (0,T]$, we have the estimates
\begin{align}
\label{stime1}
&\left|v_j(t,\bs+h)-v_j(t,\bs)-h\,\partial_s v_j(t,\bs)\right|^2\,\le\,
K_1\,\phi_2^2(t) \,|\langle y_0,e_j\rangle|^2\,h^4,\\[1mm]
\label{stime2}
&\left|\partial_s v_j(t,\bs+h)-\partial_s v_j(t,\bs)-h\,\partial^2_{ss} v_j(t,\bs)\right|^2
\\
&\nonumber\le\,
K_2\,\phi_3^2(t) \,|\langle y_0,e_j\rangle|^2\,h^4,\\[1mm]
\label{stime3}
&\left|w_j(t,\bs+h)-w_j(t,\bs)-h\,\partial_s w_j(t,\bs)\right|^2\\
&\nonumber\le\,
K_3\int_0^T\!\!\phi_2^2(t) dt\,\sup_{\theta\in (0,T)} 
|\langle f(\cdot,\theta),e_j\rangle|^2\,h^4,\\[1mm]
\label{stime4}
&\left|\partial_s w_j(t,\bs+h)-\partial_s w_j(t,\bs)-
h\,\partial^2_{ss} w_j(t,\bs)\right|^2\\
&\nonumber \le\,
K_4\int_0^T\!\!\phi_3^2(t) dt\,\sup_{\theta\in (0,T)} |\langle f(\cdot,\theta),e_j\rangle|^2\,h^4\,.
\end{align}

\noindent
From \eqref{stime1} and \eqref{stime3}, we infer that there is a constant $K_5>0$,
which is independent of $0<|h|<\bs/2$, such that
\begin{align*}
&\Bigl\|y(\bs+h)-y(\bs)-h\,\sum_{j\in\N}\partial_s y_j(\cdot,\bs)\,e_j\Bigr\|^2_{L^2(Q)}\\
&\nonumber\le\,\lim_{n\to\infty} \sum_{j=1}^n\int_0^T|y_j(t,\bs+h)-y_j(t,\bs)
-h\,\partial_s y_j(t,\bs)|^2\,dt\\
&\nonumber\le \,K_5\,\Bigl(\sum_{j\in\N}|\langle y_0,e_j\rangle|^2\,+\,\sum_{j\in\N}
f_j^2\Bigr)\,h^4\,.
\end{align*}

\noindent Hence, $\mathcal{S}$ is Fr\'echet differentiable at $\bs$ 
as a mapping from $(0,+\infty)$
into $L^2(Q)$, and the Fr\'echet derivative is given by the linear mapping
$$
h \mapsto D_s\mathcal{S}(\bs)(h)\,=\,h\,\sum_{j\in\N}\partial_s y_j(\cdot,\bs)\,e_j,
$$
as claimed. The corresponding result for the second Fr\'echet derivative follows
similarly employing the estimates \eqref{stime2} and \eqref{stime4}. This
concludes the proof of the assertion. \qed 

\section{Optimality conditions}

\noindent
In this section, we establish first-order necessary and second-order sufficient 
optimality conditions for the control problem {\bf (IP)}.  We do not address the
question of existence of optimal controls, here; this will be the subject of the forthcoming section. We have the
following result.

\vspace{3mm}
\begin{theorem}\label{CL3A}
Suppose that
that $f:\Omega\times[0,T]\to\R$ satisfies~$f(\cdot,t)\in L^2(\Omega)$, for every~$t\in[0,T]$, as well
as condition \eqref{L2-l2}. Moreover, let $y_0\in L^2(\Omega)$ be given. Then the following holds true:

\vspace{2mm}\noindent
{\bf (i)} \,\,\,If $\bs\in (0,L)$ is an optimal parameter for {\bf (IP)} and $y(\bs)$ is the associated (unique) solution to the state system \eqref{ss1}--\eqref{ss2} according to Theorem 1, then
\begin{align}
\label{necessary}
\int_0^T\!\!\int_\Omega (y(\bs)-y_Q)\,\partial_s y(\bs)\,dx\,dt \,+\,\varphi'(\bs)\,=\,0,
\end{align}
where $\,\partial_s y(\bs)$ is given by \eqref{deriS1}.

\vspace{2mm}\noindent
{\bf (ii)} \,\,If $\bs\in (0,L)$ satisfies condition \eqref{necessary} and, in addition,
\begin{align}
\label{sufficient}
\int_0^T\!\!\int_\Omega\left[
(\partial_s y(\bs))^2\,+\,(y(\bs)-y_Q)\,\partial^2_{ss} y(\bs)\right]dx\,dt
\,+\,\varphi''(\bs)\,>\,0,
\end{align}
where $\,\partial^2_{ss} y(\bs)\,$ is defined in \eqref{deriS1},
then $\bs$ is optimal for {\bf (IP)}.
\end{theorem}

\vspace{3mm}\noindent
{\sc Proof:} \,\,By Theorem \ref{II O}, the ``reduced'' cost functional 
$\,s\mapsto \mathcal{J}(s):=J(y(s),s)\,$ is twice differentiable on $(0,L)$, and it  follows directly
from the chain rule that 
\begin{align*}
\mathcal{J}'(\bs)\,&=\,\frac d{ds} \,J(y(\bs),\bs)\,=\,\partial_yJ(y(\bs),\bs)\circ D_s\mathcal{S}(\bs)
+  \partial_s J(y(\bs),\bs)\\
&=\,\int_0^T\!\!\int_\Omega (y(\bs)-y_Q)\,\partial_s y(\bs)\,dx\,dt \,+\,\varphi'(\bs)\,.
\end{align*}
Moreover,
\begin{align*}
\mathcal{J}''(\bs)\,=\,\int_0^T\!\!\int_\Omega\left[
(\partial_s y(\bs))^2\,+\,(y(\bs)-y_Q)\,
\partial^2_{ss} y(\bs)\right]\,dx\,dt
\,+\,\varphi''(\bs)\,.
\end{align*}
The assertions {\bf (i)} and {\bf (ii)} then immediately follow.\qed

\vspace{7mm}\noindent
{\bf Remark:} \,\,\,
In our framework, optimizers $\bs$ can be found by minimizing methods
(see Theorem \ref{TH 5}): in this setting,
the conditions in \eqref{phi} assure that the optimal
parameter $\bs$ lies in the open interval $(0,L)$.
Also, if $\varphi'(s)$ blows up near $0$ faster than $1/{s}$
(as it happens in the examples given in \eqref{EXA:PHI}),
solutions of \eqref{necessary} do not accumulate near $0$,
since, by \eqref{la:sec2} and \eqref{estideri1}, $$ \left|
\int_0^T\!\!\int_\Omega (y(s)-y_Q)\,\partial_s y(s)\,dx\,dt \right|
\le \| y(s)-y_Q \|_{L^2(Q)}\,\| \partial_s y(s)\|_{L^2(Q)}\le\frac{C}{s},$$
for some $C>0$.

\vspace{5mm}\noindent
{\bf Remark:} \,\,\, 
It is customary in optimal control theory to
formulate the first-order necessary optimality conditions
in terms of a variational inequality (which encodes possible
control constraints) and an adjoint state equation, while second-order
sufficient condition also involve the so-called
``$\tau$--critical cone'' (see, e.\,g., the textbook \cite{TR}). In our
situation, we can avoid these abstract concepts, since we have explicit
formulas for the relevant quantities at our disposal. Indeed, in order to
evaluate $\,y(s)\,,\,\partial_s y(s)\,,\,\partial_{ss}y(s)$, we can
use the series representations given in \eqref{sol1}
and \eqref{deriS1}. In practice, this
amounts to determining the eigenvalues $\,\lambda_j\,$ and the associated
eigenfunctions $e_j$ up to a sufficiently large index $j$,
and then to making use of the differentation formulas for the functions
\eqref{DEF:E}
for $\lambda=\lambda_j$ that are provided at the beginning
of the proof of Lemma \ref{DER:1}.
Using a standard technique (say, Newton's
method), we then can easily find an approximate minimizer of the cost
functional. Also in the case that  control constraints
$-\infty< a\le s\le b <+\infty$ are to be respected, this strategy would
still work to find interior minimizers $\bar s\in (a,b)$,
while the value of the cost at $a$ and $b$ can also be calculated.

\vspace{5mm}\noindent
{\bf Remark:} \,\,\, We recall that
in infinite dimensional setting 
conditions like \eqref{sufficient}
are not 
necessarily sufficient conditions, 
see Example 3.3 in \cite{MR3311948}. 
On the other hand, this is the case in finite dimensions.
\medskip

\begin{figure}[h]
    \centering
    \includegraphics[height=6cm]{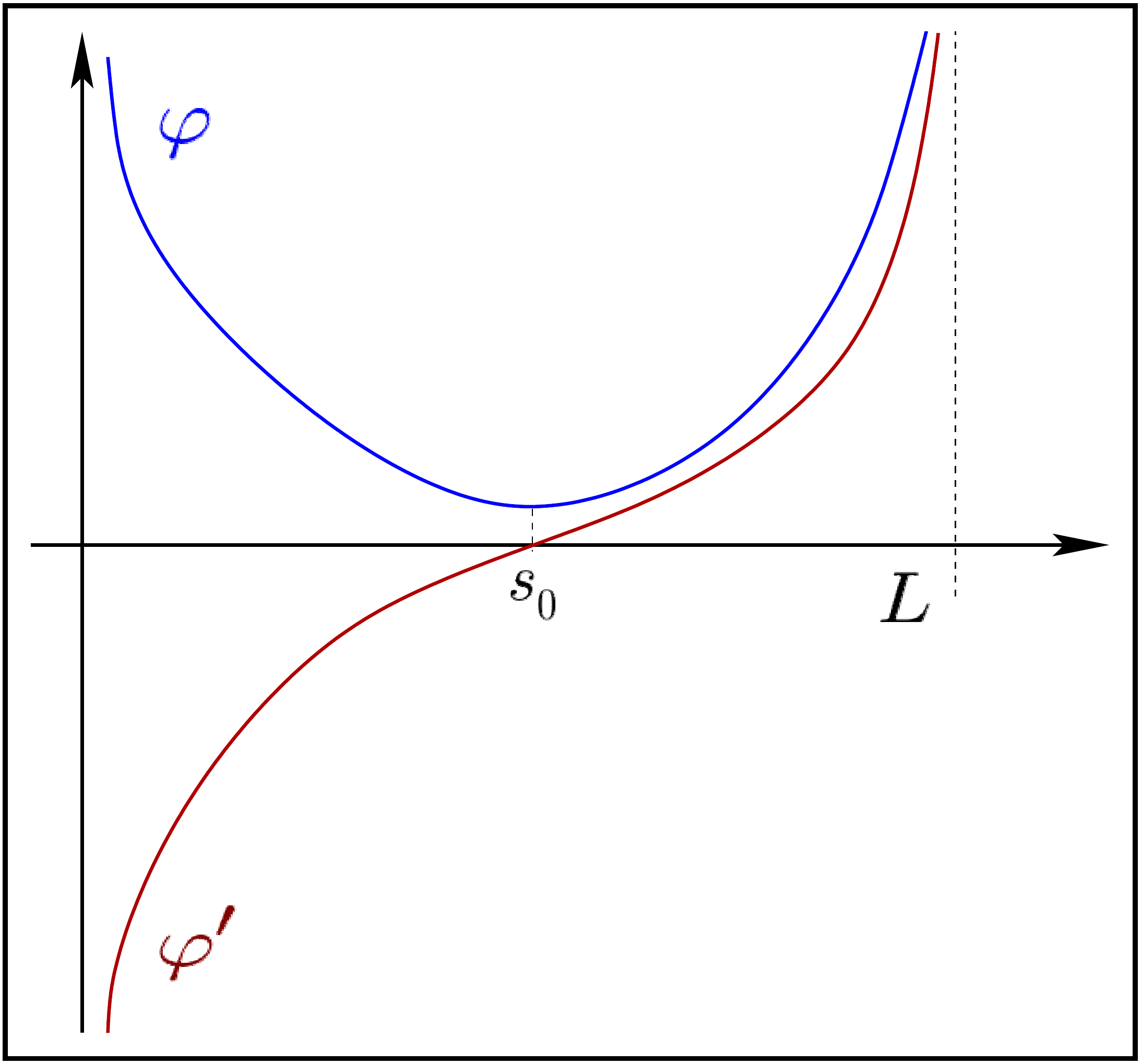}
    \caption{\small The natural cost function~$\varphi$ and its derivative.}
    \label{CORCO}
\end{figure}

To clarify Theorem~\ref{CL3A}, we now present two simple explicit examples
that outline the behavior of the optimal exponent~$\bs$
(recall~\eqref{necessary} and~\eqref{sufficient}). 
To make the arguments as simple as possible, we assume that~$\vp$ is strictly convex
and that the forcing term~$f$ is identically zero
(as a matter of fact, the functions~$\vp$ presented in~\eqref{EXA:PHI}
as examples fulfill also this convexity assumtpion). Notice that under these 
assumptions on~$\vp$ the function~$\vp$
has a unique critical point $s_0\in(0,+\infty)$, which
is a minimum (see Figure~\ref{CORCO}).
The examples are related to the fractional
Laplacian in one variable, namely, the case
of homogeneous Neumann data and the case of homogeneous Dirichlet
data on an interval. We will see that, in general, the
optimal exponent~$\bs$ differs from the minimum~$s_0$ of~$\vp$
(and, in general, it can be both larger or smaller). In a sense,
this shows that different boundary data and different target distributions~$y_Q$
influence the optimal exponent~$\bs$ and its relation with
the minimum~$s_0$ for~$\varphi$.
\medskip

\noindent{\bf Example 1.} Consider as operator~${\mathcal{L}}$
the classical~$\,-\Delta$ on the interval~$(0,\pi)$
with homogeneous Neumann data. In this case,
we can take as eigenfunctions~$e_j(x):= c_j\,\cos(j\,x)$,
where~$c_j\in\R\setminus\{0\}$ is a normalizing constant,
and~$j=0,1,2,3,\dots$. The eigenvalue corresponding to~$e_j$
is~$\lambda_j=j^2$.

Now let, with a fixed~$j_0\in\N$, where $j_0>1$, and~$\epsilon\in\R$, 
\begin{equation*}
y_0(x):= 1 + \epsilon\, e_{j_0}(x) \quad\forall\,x\in [0,\pi].
\end{equation*}
Then it is easily verified that for every $s>0$ the 
unique solution to \eqref{ss1}, \eqref{ss2} is given by
$$
y(s)(x,t) = 1 + \epsilon \, e_{j_0}(x)\, e^{-j_0^{2s}\,t} \quad\forall\,
(x,t)\in\overline{Q}.
$$
We now make the special choice $y_Q(x,t):=1$ for the target function.
We then observe that
$$ \partial_s y(s)(x,t) = -2\epsilon \,j_0^{2s}\, \ln(j_0)\,
t\, e_{j_0}(x)\, e^{-j_0^{2s}\,t},$$
and therefore, using the substitution~$\vartheta:=j_0^{2s}\,t$,
\begin{eqnarray*}
&& \int_0^T\!\!\int_\Omega (y(s)-y_Q)\,\partial_s y
(s) \,dx\,dt \\
&=& 
-2\epsilon^2 \,j_0^{2s}\, \ln(j_0)\,
\int_0^T\!\!\int_\Omega 
t\, e_{j_0}^2(x)\, e^{-2j_0^{2s}\,t}
\,dx\,dt \\
&=& 
-2\epsilon^2 \,j_0^{2s}\, \ln(j_0)\,
\int_0^T t\, e^{-j_0^{2s}\,t}\,dt \\
&=&
-2\epsilon^2 \,j_0^{-2s}
\,\ln(j_0)\,
\int_0^{j_0^{2s}T} \vartheta\, e^{-2\vartheta}\,dt.
\end{eqnarray*}
As a consequence, condition~\eqref{necessary} becomes,
in this case,
\begin{equation}\label{CA-1}
\varphi'(\bs) =
2\epsilon^2 \,j_0^{-2\bs}
\,\ln(j_0)
\int_0^{j_0^{2\bs}T} \vartheta\, e^{-2\vartheta}\,dt.\end{equation}
If~$\epsilon=0$ (and when~$j_0\to+\infty$), then the identity in~\eqref{CA-1}
reduces to~$\varphi'(\bs) =0$; that is, in this case
the ``natural'' optimal exponent~$s_0$ coincides with
the optimal exponent~$\bs$ given by the full cost functional
(that is, in this case the external conditions given by
the exterior forcing term and the resources do not alter the
natural diffusive inclination of the population).

\begin{figure}
    \centering
    \includegraphics[height=6cm]{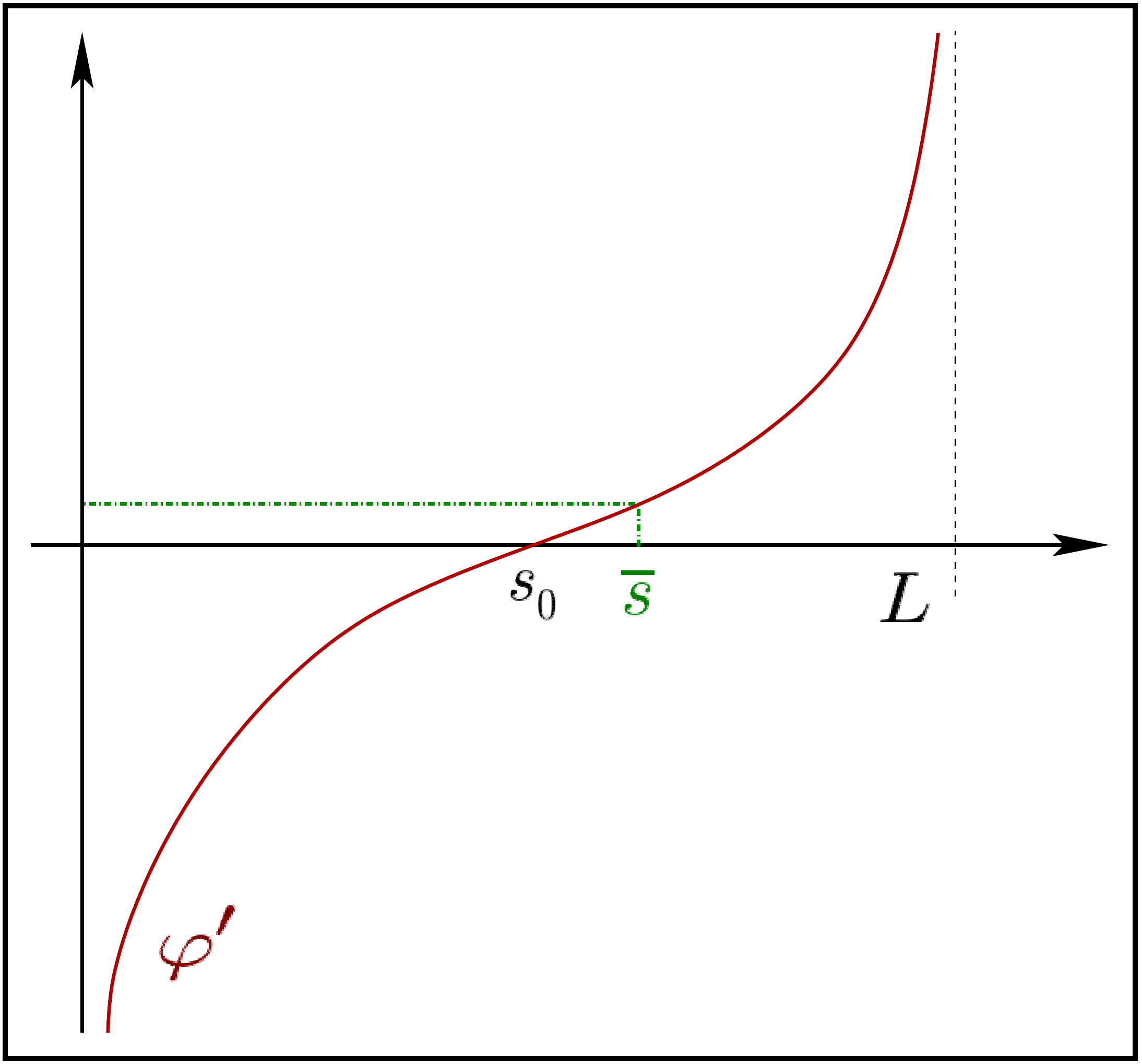}
    \caption{\small The optimal exponent~$\bs$ in Example~1.}
    \label{FIG:NEUMANN}
\end{figure}

But, in general, for fixed $\epsilon\ne0$ and $j_0>1$, 
the identity in~\eqref{CA-1} gives that~$\varphi'(\bs) >0$.
This, given the convexity of~$\varphi$, implies that~$\bs >s_0$,
i.\,e., the optimal exponent given by the cost functional 
is larger than the natural one (see Figure~\ref{FIG:NEUMANN}).
\medskip

\noindent{\bf Example 2.} Now we consider as operator~${\mathcal{L}}$
the classical~$-\Delta$ on the interval~$(0,\pi)$ 
with homogeneous Dirichlet data.
In this case,
we can take as eigenfunctions~$e_j(x):= c_j\,\sin(j\,x)$,
where~$c_j\in\R\setminus\{0\}$ is a normalizing constant,
and~$j=1,2,3,\dots$. The eigenvalue corresponding to~$e_j$
is~$\lambda_j=j^2$.

\begin{figure}[h]
    \centering
    \includegraphics[height=6cm]{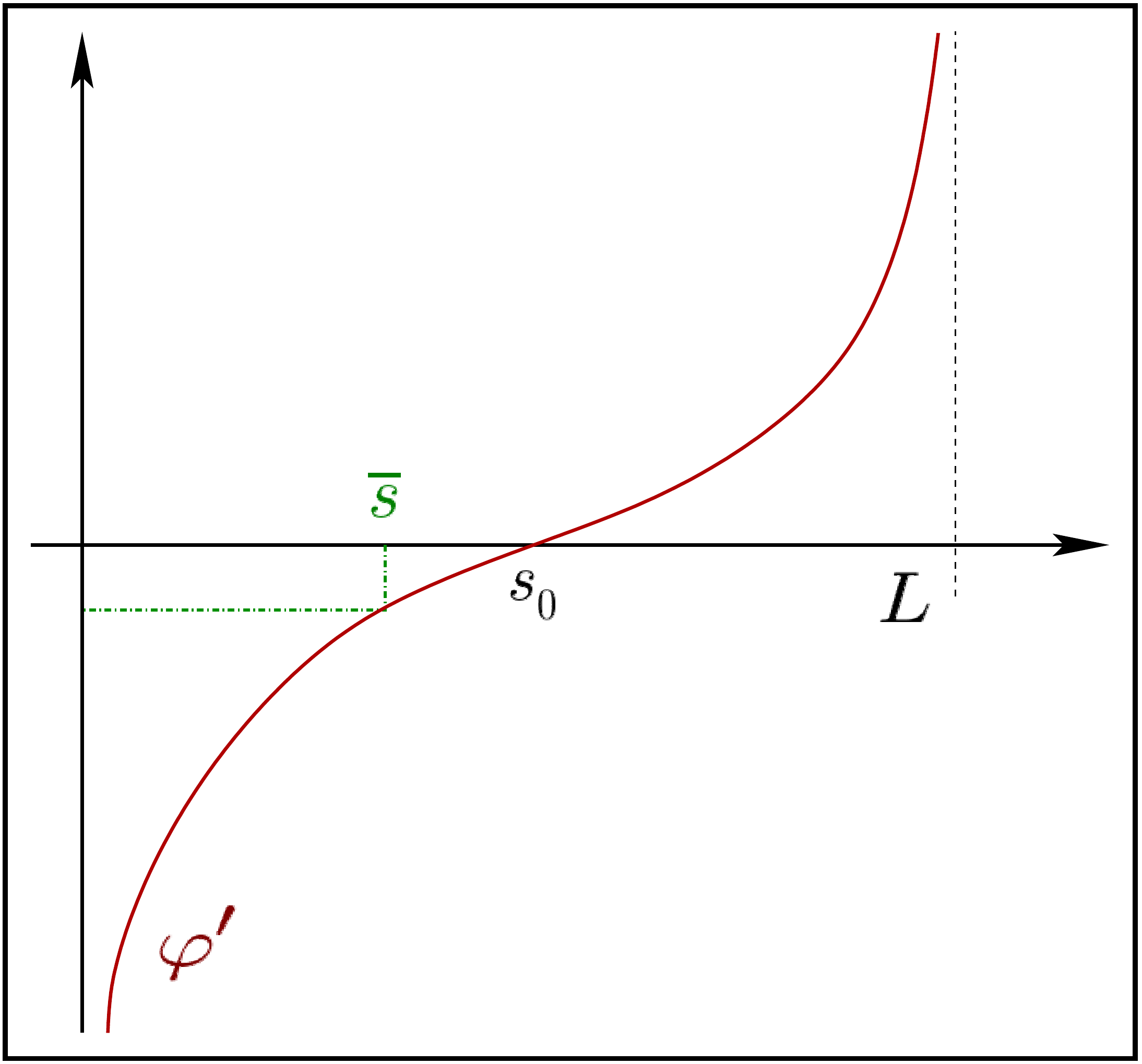}
    \caption{The optimal exponent~$\bs$ in Example~2.}
    \label{FIG:DIRI}
\end{figure}

\noindent
For fixed~$j_0\in\N$ with~$j_0\ge 1$, and~$\epsilon\in\R$, we set
$$
y_0(x) := \epsilon\, e_{j_0}(x)\quad\forall\,x\in [0,\pi]\,.
$$
Then, for every $s>0$, the corresponding solution is given by
$$
y(s)(x,t) = \epsilon \,e_{j_0}(x)\, e^{-j_0^{2s}\,t}
\quad\forall\,(x,t)\in\overline{Q}\,.
$$
Now, let $y_Q(x,t):=\epsilon\, e_{j_0}(x)$ for $(x,t)\in Q$. We have
$$ \partial_s y(s)(x,t) = -2\epsilon \,j_0^{2s}\, \ln(j_0)\,
t\, e_{j_0}(x)\, e^{-j_0^{2s}\,t},$$
and therefore, using the substitution~$\vartheta:=j_0^{2s}\,t$,
\begin{eqnarray*}
&& \int_0^T\!\!\int_\Omega (y(s)-y_Q)\,\partial_s y
(s) \,dx\,dt \\
&=& -2\epsilon^2\,j_0^{2s}\, \ln(j_0)\, \int_0^T\!\!\int_\Omega
t\,e_{j_0}^2(x)\, \big( e^{-j_0^{2s}\,t} -1\big)\,e^{-j_0^{2s}\,t}
\,dx\,dt  \\ &=&
-2\epsilon^2\,j_0^{2s}\, \ln(j_0)\, \int_0^T
t\, \big( e^{-j_0^{2s}\,t} -1\big)\,e^{-j_0^{2s}\,t}
\,dt \\ &=&
-2\epsilon^2\,j_0^{-2s}\,\ln(j_0)\, \int_0^{j_0^{2s} T}
\vartheta\, \big( e^{-\vartheta} -1\big)\,e^{-\vartheta}
\,d\vartheta.
\end{eqnarray*}
So, in this case,
condition~\eqref{necessary} becomes
\begin{equation}\label{CA-2}
\varphi'(\bs) =
2\epsilon^2\,j_0^{-2\bs}\,\ln(j_0)\, \int_0^{j_0^{2\bs} T}
\vartheta\, \big( e^{-\vartheta} -1\big)\,e^{-\vartheta}
\,d\vartheta.
\end{equation}
If~$\epsilon=0$ (and when~$j_0\to+\infty$), then the identity in~\eqref{CA-2}
reduces to~$\varphi'(\bs) =0$, which boils down to~$\bs=s_0$.
But if~$\epsilon\ne0$ and~$j_0\ge1$, then 
the identity in~\eqref{CA-2} gives that~$\varphi'(\bs) <0$.
By the convexity of~$\varphi$, this implies that~$\bs <s_0$,
i.\,e., the optimal exponent given by the full cost functional 
is in this case smaller  than the natural one (see Figure~\ref{FIG:DIRI}).
\medskip

We observe that, in the framework
of Examples~1 and~2, the effect
of a larger~$s$ is to ``cancel faster'' the higher order
harmonics in the solution~$y$; since these harmonics
are related to ``wilder oscillations'', one may think that
the higher $s$ becomes, the bigger the smoothing effect 
is. In this regard,
roughly speaking, a larger~$s$ ``matches better'' with a constant
target function~$y_Q$  and a smaller~$s$ with an oscillating one
(compare again Figures~\ref{FIG:NEUMANN} and~\ref{FIG:DIRI}).

We also remark that when~$j_0\ge2$ in Example~2 (or if~$\epsilon$
is large in Example~1), the solution~$y$ is not positive.
On the one hand, this seems to reduce the problem, in this case,
to a purely mathematical question, since if~$y$ represents
the density of a biological population, the assumption~$y\ge0$
seems to be a natural one. On the other hand, there are other models
in applied mathematics in which the condition~$y\ge0$ is not assumed:
for instance, if~$y$ represents the availability of specialized
workforce in a given field, the fact that~$y$ becomes negative
(in some regions of space, at some time) translates into the fact
that there is a lack of this specialized 
workforce (and, for example, non-specialized workers have to
be used to compensate this lack).

The use of mathematical models to deal with problems in the job market
is indeed an important topic of contemporary research, see e.g.
\cite{Stew} and the references therein.\medskip

The models arising in the (short time) job market also provide
natural examples in which the birth/death effects in
the diffusion equations are negligible.
\vspace{3mm}

\section{Existence and a compactness lemma} 
\noindent
In this section, we establish an existence result for the identification problem {\bf (IP)}. We make the
following general assumption for the initial datum $y_0$:
\begin{equation}\label{initial} 
\sup_{s\in(0,L)} \|y_0\|_{{\mathcal{H}}^s}<+\infty.
\end{equation}

\noindent
{\bf Remark:} \,\,\,We remark that the condition \eqref{initial} can be very restrictive if $L$ is large. Indeed,
we obviously have $\,\lambda_j^{2s}\le 1\,$ for $\,\lambda_j\le 1$, and for $\,\lambda_j>1\,$ the
function $\,\,s\mapsto \lambda_j^{2s}\,$ is strictly increasing. From this it follows that
\eqref{initial} is certainly fulfilled for a finite $L$ if only $\|y_0\|_{{\mathcal{H}}^L}<+\infty$, that is, if 
$y_0\in {\mathcal{H}}^L$. 

For an example, consider the prototypical case when ${\mathcal{L}}=-\Delta$ with zero
Dirichlet boundary condition. Then the choice $L=\frac 12$ leads to the
requirement $y_0\in H^1_0(\Omega)$, while
for the choice $L=1$ we must have $y_0\in H^2(\Omega)\cap H^1_0(\Omega)$: indeed, if $\,\{\lambda_j\}_{j\in\N}\,$
are the corresponding eigenvalues with associated orthogonal eigenfunctions $\{e_j\}_{j\in\N}$,
normalized by $\|e_j\|_{L^2(\Omega)}=1$\, for all $j\in\N$, then it is readily
verified that the set  $\{\lambda_j^{-1/2}e_j\}_{j\in\N}$ forms an orthonormal basis in the Hilbert
space $\bigl(H^1_0(\Omega), \langle \cdot,\cdot\rangle_1\bigr)$ with respect to the inner product
$\,\langle u,v\rangle_1:=\xinto\nabla u\cdot\nabla v\,dx$. Therefore, if $y_0\in H^1_0(\Omega)$, it follows 
from Parseval's identity and integration by parts that
\begin{align*}
+\infty\,&>\,\|y_0\|^2_{H^1_0(\Omega)}\,=\,\sum_{j\in\N}
\left|\langle y_0,\lambda_j^{-1/2}\,e_j\rangle_1 \right|^2\\
&=\,\sum_{j\in\N}\frac 1{\lambda_j}\Bigl|\xinto\nabla y_0\cdot\nabla e_j\,dx\Bigr|^2
\,=\,\sum_{j\in\N}\frac 1{\lambda_j}\Bigl|-\xinto y_0\,\Delta e_j\,dx \Bigr|^2\\
&=\,\sum_{j\in\N}\lambda_j\,|\langle y_0,e_j\rangle|^2\,=\,\|y_0\|^2_{\mathcal{H}^{1/2}}\,.
\end{align*}
The case $L=1$ is handled similarly. It ought to be clear that with increasing $L$ the 
condition \eqref{initial} imposes ever higher regularity postulates on $y_0$. On the other hand,
\eqref{initial} is obviously satisfied for every finite $L>0$ if $y_0$ belongs to the set
of finite linear combinations of the eigenfunctions $\{e_j\}_{j\in\N}$, that is, on a 
dense subset of $L^2(\Omega)$.

\vspace{6mm}
We now give sufficient conditions that guarantee the existence of a solution to the 
optimal control problem {\bf (IP)}.

\vspace*{3mm}
\begin{theorem}\label{TH 5}
Suppose
that $f:\Omega\times[0,T]\to\R$ satisfies~$f(\cdot,t)\in L^2(\Omega)$, for every~$t\in[0,T]$, as well
as condition \eqref{L2-l2}. Moreover, let $y_0\in L^2(\Omega)$ satisfy the condition \eqref{initial}.
If $\,\lambda_j\nearrow +\infty\,$ as $\,j\to +\infty$, then the control problem {\bf (IP)} has a 
solution, that is, $\mathcal{J}$\,attains a
minimum in $(0,+\infty)$.  
\end{theorem}

\vspace{5mm}
Before proving the existence result, we establish an auxiliary compactness lemma, which is of some
interest in itself, since it acts between spaces with different fractional 
coefficients $s$. 

\begin{lemma}\label{compactness}
Assume that the sequence $\{\lambda_k\}_{k\in\N}$ of eigenvalues of $\,\mathcal{L}\,$ satisfies
$\,\lambda_k\nearrow +\infty$ as~$k\to\infty$, and assume
that the sequence $\,\{s_k\}_{k\in\N}\subset (0,+\infty)$ satisfies $\,s_k\to\bs\,$
as $k\to\infty$, for some~$\,\bs\in(0,+\infty)\cup\{+\infty\}$.
Moreover, let a sequence $\,\{y_k\}_{k\in\N}\,$ be given such that 
$\,y_k\in L^2(0,T;{\mathcal{H}^{s_k}})\,$ and 
$\,\partial_t y_k\in L^2(Q)$, for all $k\in\N$, as well as
\begin{align}\label{CST} 
& \sup_{k\in\N}\Big( 
\| y_k \|_{L^2(Q)}\,+\, 
\| y_k \|_{L^2(0,T;{\mathcal{H}^{s_k}})}\Big)<+\infty, \quad
\mbox{and }\;\, \\[1mm]
&\nonumber
\sup_{k\in\N}\| \partial_t y_k\|_{L^2(Q)}<+\infty .
\end{align}
Then~$\,\{y_k\}_{k\in\N}\,$ contains a subsequence that converges strongly 
in~$\,L^2(Q)$.
\end{lemma}

\vspace{4mm}\noindent
{\sc Proof:} \,\,For fixed~$j\in\N$, we define
$$ y_{k,j}(t):=\int_\Omega y_k(x,t)\,e_j(x)\,dx.$$
Notice that
\begin{eqnarray*} 
&& \int_0^T |\partial_t y_{k,j}(t)|^2\,dt
\,\le\, \int_0^T\!\!\left(
\int_\Omega |\partial_t y_k(x,t)|\,|e_j(x)|\,dx\right)^2 dt\\
&&\quad\le 
\int_0^T\!\!\left(
\int_\Omega |\partial_t y_k(x,t)|^2\,dx\right)dt \,=\,
\| \partial_t y_k\|_{L^2(Q)}^2,
\end{eqnarray*}
which is bounded uniformly in~$k$, thanks to~\eqref{CST}. 

Hence, we obtain a
bound in~$H^1(0,T)$ for~$y_{k,j}$, which is uniform  with respect to~$k\in\N$, for every $j \in\N$.
Owing to the compactness of the embedding $H^1(0,T)\subset C^{1/4}([0,T])$,
the sequence~$\{y_{k,j}\}_{k\in\N}$ thus forms for every $j\in\N$  a compact subset~$C_j$
of~$C^{1/4}([0,T])$.

Therefore, the infinite string~$\big(\{y_{k,1}\}_{k\in\N},\,
\{y_{k,2}\}_{k\in\N},\dots\big)$ lies in~$C_1\times C_2\times\dots$,
which, by virtue of Tikhonov's Theorem,  is compact in the product space
$$C^{1/4}([0,T])\times C^{1/4}([0,T])\times\dots\,.$$
Hence, there is a subsequence
(denoted by the index~$k_m$), which converges in
this product space to an infinite string of the form~$\big(y^*_{1},y^*_2,\dots\big)$.
More explicitly, we have that~$y^*_j\in C^{1/4}([0,T])$, for any~$j\in\N$, and
\begin{equation}\label{0s72uuJ}
\lim_{m\to\infty} \| y_{k_m,j} -y^*_j\|_{C^{1/4}([0,T])}=0 \quad\mbox{for every $j\in\N$.}
\end{equation}

\vspace{2mm}\noindent
We then define
$$ y^*(x,t):=\sum_{j\in\N} y^*_j(t)\,e_j(x)$$
and  claim that
\begin{equation}\label{78HHA}
y_{k_m}\to y^*\quad\mbox{strongly  in}\,L^2(Q).
\end{equation}
To prove this claim,
we fix~$\epsilon\in(0,1)$ and choose~$j_*\in\N$ so large  that
\begin{equation}\label{9wkUUa}
{\mbox{$ \lambda_j\ge\epsilon^{-1}$ for any~$j\ge j_*$.}}\end{equation}
Then, by~\eqref{0s72uuJ}, we may also fix~$m_*\in\N$ large enough, so that
for any~$m\ge m_*$ it holds that
$$ s_{k_m} \ge \min\left\{ 1,\,\frac{\bs}{2}\right\}=:\sigma,$$
as well as
$$ \| y_{k_m,j} -y^*_j\|_{C^{1/4}([0,T])}\le \frac{\epsilon}{j_*+1}\qquad
{\mbox{ for every }}j<j_*.
$$

\vspace{2mm}\noindent
Now, let $t\in(0,T)$ be fixed. Then, for any~$m\ge m_*$,
\begin{align}\label{78:8YUU12}
& \| y^*(\cdot, t)-y_{k_m}(\cdot,t)\|^2_{L^2(\Omega)}\,=\,
\sum_{j\in\N} |y^*_j(t) -y_{k_m,j}(t)|^2\\
&\nonumber\le\, \sum_{{j\in\N}\atop{j<j_*}} |y^*_j(t) -y_{k_m,j}(t)|^2
\,+\,4\sum_{{j\in\N}\atop{j\ge j_*}}\big( |y^*_j(t)|^2 +|y_{k_m,j}(t)|^2\big)\\
&\nonumber\le\, \epsilon
+4\sum_{{j\in\N}\atop{j\ge j_*}}\big( |y^*_j(t)|^2 +|y_{k_m,j}(t)|^2\big).
\end{align}
Moreover, by~\eqref{9wkUUa},
for any~$\ell\in\N$,
\begin{align}\label{UI:0011}
& \sum_{ {j\in\N}\atop{j_*\le j\le j_*+\ell}} |y_{k_m,j}(t)|^2
\,\le\, \sum_{{j\in\N}\atop{j_*\le j\le j_*+\ell}} \epsilon^{2s_{k_m}}
\lambda_j^{2s_{k_m}}|y_{k_m,j}(t)|^2\\[2mm]
&\nonumber\le\, \epsilon^{2\sigma}\, \|y_{k_m}
\|_{{\mathcal{H}}^{s_{k_m}}}^2\,\le\, \epsilon^{2\sigma}M,
\end{align}
for some~$M>0$, where the last inequality follows from~\eqref{CST}.
Hence, by virtue of \eqref{0s72uuJ}, taking limit as $\,m\to\infty$, we obtain that
\begin{equation}\label{UI:0012}
\sum_{ {j\in\N}\atop{j_*\le j\le j_*+\ell}} |y^*_{j}(t)|^2
\le \epsilon^{2\sigma}M
.\end{equation}
Therefore, letting~$\ell\to\infty$ in~\eqref{UI:0011}
and~\eqref{UI:0012}, we find that
$$ \sum_{ {j\in\N}\atop{j\ge j_*}} |y_{k_m,j}(t)|^2
\le \epsilon^{2\sigma}M\;\;{\mbox{ and }}\;\;
\sum_{ {j\in\N}\atop{j\ge j_*}} |y^*_{j}(t)|^2
\le \epsilon^{2\sigma}M.$$

\vspace{2mm}\noindent
Insertion of  these bounds in~\eqref{78:8YUU12} then
yields that
$$  \| y^*(\cdot, t)-y_{k_m}(\cdot,t)\|_{L^2(\Omega)}^2
\le \epsilon
+8\epsilon^{2\sigma}M,$$
as long as~$m\ge m_*$.
By taking $\epsilon$ arbitrarily small, we conclude the
validity of~\eqref{78HHA}
and thus of the assertion of the lemma.
\qed

\vspace*{7mm} \noindent
{\sc Proof of Theorem \ref{TH 5}:} \,\,The proof is a combination of the Direct Method
with the regularity results proved in Theorem~1 
and the compactness argument stated in Lemma~6.
First of all, we observe that $\mathcal{J}(\frac L2)<+\infty$  if $0<L<+\infty$,
while ${\mathcal{J}}( \frac 12)<+\infty$ if $L=+\infty$. Hence, owing to \eqref{phi}, we have  
$$
0<\inf_{0<s< L} \mathcal{J}(s)<+\infty.
$$

\vspace{2mm}\noindent
Now, we pick a  minimizing sequence $\{s_k\}_{k\in\N}\subset (0,L)$
and consider, for every $k\in\N$, the (unique) solution $\,y_k:=\mathcal{S}(s_k)=y(s_k)$
to the state system \eqref{ss1}, \eqref{ss2} associated with $s=s_k$. 
We may without loss of generality assume that
$$ {\mathcal{J}}(s_k)\le 1+{\mathcal{J}}(s^*) \quad\forall \,k\in\N , 
$$
where \,$s^*:=\frac L2$\, if $L<+\infty$ and $\,s^*:=\frac 12$\, otherwise. We then infer that
\begin{equation}
\label{b1}
\| y_k\|_{L^2(Q)}\,+\,\varphi(s_k)\,\le\,C_1 \quad\forall\,k\in\N,
\end{equation}
where, here and in the following, we denote by $C_i$, $i\in\N$, constants that may depend on the
data of the problem but not on $k$. 
In particular, by~\eqref{phi}, the sequence $\,\{s_k\}_{k\in\N}\,$ is bounded, and we may without loss
of generality assume that \,$s_k\to \bs$\, for some~$\,\bs\in(0,L)$.

\noindent
Also, by virtue of \eqref{la:sec2} and \eqref{initial}, we obtain that
\begin{equation}\label{U8:TY:A}
\| \partial_t y_k\|_{L^2(Q)}\,+\,
\| y_k\|_{L^2(0,T;{\mathcal{H}}^s)}\,\le\,C_2\,,
\end{equation}
whence,  in particular,
\begin{equation}
\label{b2}
\sum_{j\in\N} \int_0^T |\langle \partial_t y_k(\cdot,t), e_j\rangle|^2\,dt\,\le\,C_3\quad\forall\,k\in\N.
\end{equation}

\noindent
Thus, using the compactness result of Lemma 6, we can select a subsequence, which is again indexed by $k$,
such that there is some $\by\in H^1(0,T;L^2(\Omega))$ satisfying
\begin{align}
y_k&\to \by \quad\mbox{strongly in $L^2(Q)$ and pointwise a.\,e. in $Q$}, \\
\nonumber y_k&\to \by \quad\mbox{weakly in $H^1(0,T;L^2(\Omega))$}.
\end{align}
Therefore, we can infer from \eqref{b2} that
\begin{equation}
\label{b3}
\sum_{j\in\N} \int_0^T |\langle \partial_t \by(\cdot,t), e_j\rangle|^2\,dt\,\le\, C_3.
\end{equation}
We now claim that~$\by=y(\bs)$, that is, that $\by$  is the (unique) solution to the state system
associated with $s=\bs$. To this end, it suffices to show that $\by$ satisfies the
conditions \eqref{OV}--\eqref{SL3}, since then the claim follows exactly in the same  way as
uniqueness was established in the proof of Theorem 1; in this connection, observe that
for this argument the validity of \eqref{SL1} was not needed.

To begin with, we fix $j\in\N$. We conclude from \eqref{b2} that it holds 
that
$$ \int_0^T |\partial_t \langle y_k(\cdot,t),e_j\rangle|^2\,dt\,\le\,C_4\quad\forall \,k\in\N\,.
$$
Hence,  the sequence formed by the mappings $\,t\mapsto\langle y_k(\cdot,t),e_j\rangle$ 
is  a bounded subset of  $H^1(0,T)$. Hence, its weak limit, which is given by the mapping
$\,t\mapsto \langle \by(\cdot,t),e_j\rangle$, belongs to $H^1(0,T)$ and is thus
absolutely continuous, which implies that \eqref{SL2} holds true for $\by$.

Moreover, by virtue of the continuity of the embedding $H^1(0,T)\subset C^{1/2}([0,T])$,
we can infer from the Arzel\`{a}--Ascoli
Theorem that the  convergence of the sequence 
$\,\{\langle y_k(\cdot,t),e_j\rangle\}_{k\in\N}\,$
is uniform on $[0,T]$. Therefore, 
to any fixed~$\epsilon>0$ there exists some $\,k_\epsilon\in\N\,$ such that,
for~$\,k\ge k_\epsilon$,
\begin{eqnarray*}
&& \big| \langle \by(\cdot,t),e_j\rangle
-\langle y_0,e_j\rangle
\big|\\
&&\le\, 
\big| \langle \by(\cdot,t),e_j\rangle
-\langle y_k(\cdot,t),e_j\rangle
\big|
\,+\, \big| \langle y_k(\cdot,t),e_j\rangle
-\langle y_0,e_j\rangle
\big|\\
&&\le\, \big| \langle y_k(\cdot,t),e_j\rangle
-\langle y_0,e_j\rangle
\big| \,+\,\epsilon.
\end{eqnarray*}
Hence, taking the limit in~$t$, and then letting $\epsilon\searrow0$,
we obtain that $\by$ fulfills \eqref{OV}.

Now we use the fact that the mapping \,$t\mapsto
\langle y_k(\cdot,t),e_j\rangle\,$ belongs to~$H^1(0,T)$ to write~\eqref{SL3}
in the weak sense. We have, for any test function~$\Psi\in C^\infty_0(0,T)$, 
\begin{eqnarray*}
&& -\int_0^T \langle y_k(\cdot,t),\,e_j\rangle\,\partial_t\Psi(t)\,dt
+\lambda_j^{s_k} \int_0^T\langle y_k(\cdot,t),\,e_j\rangle\,\Psi(t)\,dt
\\&&\quad=\int_0^T\langle f(\cdot,t),\,e_j\rangle\,\Psi(t)\,dt.
\end{eqnarray*}
Passage to the limit as $k\to\infty$ then yields the identity
\begin{eqnarray*}
&& -\int_0^T \langle \by(\cdot,t),\,e_j\rangle\,\partial_t\Psi(t)\,dt
+\lambda_j^{\bs} \int_0^T\langle \by(\cdot,t),\,e_j\rangle\,\Psi(t)\,dt
\\&&\quad=\int_0^T\langle f(\cdot,t),\,e_j\rangle\,\Psi(t)\,dt.
\end{eqnarray*}
This, and the fact that the mapping $\,t\mapsto
\langle \by(\cdot,t),e_j\rangle\,$ belongs to the space \,$H^1(0,T)$,
 give~\eqref{SL3} (recall, for instance, Theorem~6.5
in~\cite{lieb}).

In conclusion, it holds $\by=y(\bs)$, and thus the pair $(\bs,\by)$
is admissible for the problem {\bf (IP)}. By the weak sequential semicontinuity of the cost functional,
$\bs$ is a minimizer of $\mathcal{J}$. This concludes the proof of the
assertion. \qed

\section*{Acknowledgments}
\noindent
EV was supported by
ERC grant 277749 ``EPSILON Elliptic Pde's and Symmetry of Interfaces and
Layers for Odd Nonlinearities'' and 
PRIN grant
201274FYK7 ``Critical Point Theory
and Perturbative Methods for Nonlinear Differential Equations''.

\bibliographystyle{siamplain}

\begin{thebibliography}{99}
\bibitem{AO1}
{\sc H. Antil, E. Ot\'arola},
{\em A FEM for an optimal control problem of fractional powers of elliptic
operators}. Preprint arXiv:1406.7460v3 [math.OC] 18 April 2015.

\bibitem{AO2}
{\sc H. Antil, E. Ot\'arola, A. J. Salgado},
{\em A fractional space-time optimal control problem: analysis and discretization}. 
Preprint arXiv:1504.00063v1 [math.OC] 31 March 2015.

\bibitem{AO3}
{\sc H. Antil, E. Ot\'arola, A. J. Salgado},
{\em Some applications of weighted norm inequalities to the analysis of optimal control
problems}. Preprint arXiv:1505.03919v1 [math.OC] 14 May 2015.

\bibitem{Bors}
{\sc D. Bors}, {\em Optimal control of nonlinear systems governed by Dirichlet fractional Laplacian
in the minimax framework}. Preprint arXiv:1509.01283v1 [math.AP] 3 Sep 2015.

\bibitem{brezis}
{\sc H. Brezis},
{\em Functional analysis, Sobolev spaces and partial differential equations}. 
Universitext. Springer, New York (2011).

\bibitem{RO}
{\sc R. S. Cantrell, C. Cosner, Y. Lou},
{\em Advection-mediated coexistence of competing species}.
Proc. Roy. Soc. Edinburgh Sect. A 137, no. 3, 497-518 (2007).

%\juerg{\bibitem{CL}
%{\sc E. Casas, P. I. Kogut, G. Leugering}, {\em Approximation of optimal control problems in the
%coefficient for the $p$-Laplace equation. I. Convergence results.} To appear in SIAM J. Control. Optim.
%}

\bibitem{MR3311948}
{\sc E. Casas, F. Tr{\"o}ltzsch},
{\em Second order optimality conditions and their role in {PDE}
control}.
Jahresber. Dtsch. Math.-Ver. 117, no. 1, 3-44 (2015).

\bibitem{F12}
{\sc A. Friedman},
{\em PDE problems arising in mathematical biology}.
Netw. Heterog. Media 7, no. 4, 691-703 (2012).

\bibitem{H10}
{\sc N. E. Humphries, N. Queiroz, J. R. M. Dyer,
N. G. Pade, M. K. Musyl, K. M. Schaefer, D. W. Fuller,
J. M. Brunnschweiler,
T. K. Doyle, J. D. R. Houghton,	G. C. Hays,
C. S. Jones, L. R. Noble, V. J. Wearmouth,
E. J. Southall, D. W. Sims},
{\em Environmental context explains L\'evy
and Brownian movement patterns of
marine predators}.
Nature 465, 1066-1069 (2010).

\bibitem{lieb}
{\sc E. H. Lieb, M. Loss},
{\em Analysis}.
Second edition. Graduate Studies in Mathematics. 
American Mathematical Society, Providence (2001).

\bibitem{MV}
{\sc A. Massaccesi, E. Valdinoci},
{\em Is a nonlocal diffusion stra\-tegy 
convenient for biological populations in competition?}
{\tt http://arxiv.org/pdf/1503.01629.pdf}

\bibitem{P13}
{\sc E. Montefusco, B. Pellacci, G. Verzini},
{\em Fractional diffusion with Neumann boundary conditions: the logistic equation}. 
Discrete Contin. Dyn. Syst. Ser. B 18, no. 8, 2175-2202 (2013).

\bibitem{KO}
{\sc J. J. Koliha}, 
{\em A fundamental theorem of calculus for Lebesgue integration}. 
Am. Math. Mon. 113, no. 6, 551-555 (2006).

\bibitem{NAZA3}
{\sc R. Musina, A. I. Nazarov},
{\em On fractional Laplacians - 3}. To appear on 
ESAIM Control Optim. Calc. Var.

\bibitem{TWO}
{\sc R. Servadei, E. Valdinoci},
{\em On the spectrum of two different fractional operators}.
Proc. Roy. Soc. Edinburgh Sect. A 144, no. 4, 831-855 (2014).

\bibitem{Stew}
{\sc B. D. Stewart, D. B. Webster, S. Ahmad, J. O. Matson},
{\em Mathematical models for developing a flexible workforce},
Intern. Journ. of Production Economics 36, no. 3, 243-254 (1994).

\bibitem{TR}
{\sc F. Tr\"oltzsch}, {\em Optimal Control of Partial Differential
Equations: Theory, Methods and Applications}. Graduate Studies in
Mathematics Vol. 112. American Mathematical Society, Providence,
Rhode Island (2010).

\bibitem{W96}
{\sc G. M. Viswanathan, V. Afanasyev, S. V. Buldyrev, 
E. J. Murphy, P. A. Prince, H. E. Stanley},
{\em L\'evy flight search patterns of wandering albatrosses}.
Nature 381, 413-415 (1996).

\end{thebibliography}

\end{document}